\documentclass[envcountsame,envcountsect]{svmult}

\usepackage[utf8]{inputenc}
\usepackage{mathptmx}       
\usepackage{helvet}         
\usepackage{courier}        

\DeclareMathAlphabet{\mathcal}{OMS}{cmsy}{m}{n} 
\usepackage{graphicx}
\usepackage{tikz}
\usetikzlibrary{cd}
\usepackage{subcaption}
\usepackage{amsmath,amssymb,mathrsfs}
\usepackage{xcolor}

\usepackage[bookmarks=false]{hyperref}

\DeclareMathOperator{\PSL}{PSL}
\DeclareMathOperator{\SL}{SL}
\DeclareMathOperator{\Isom}{Isom}

\DeclareMathOperator{\tr}{tr}
\DeclareMathOperator{\QH}{QH}

\DeclareMathOperator{\Teich}{Teich}

\DeclareMathOperator{\MCG}{Mod}
\DeclareMathOperator{\Mod}{\mc{M}}
\DeclareMathOperator{\interior}{int}
\DeclareMathOperator{\ML}{\mc{ML}}
\DeclareMathOperator{\PML}{\mc{PML}}

\newcommand{\mc}{\mathcal}
\newcommand{\CC}{\mathbb{C}}
\newcommand{\HH}{\mathbb{H}}
\newcommand{\ZZ}{\mathbb{Z}}
\newcommand{\RR}{\mathbb{R}}
\newcommand{\QQ}{\mathbb{Q}}
\newcommand{\NN}{\mathbb{N}}
\newcommand{\MM}{\mathbb{M}}
\newcommand{\PP}{\mathbb{P}}
\newcommand{\BB}{\mathbb{B}}
\renewcommand\rightmark{Kleinian groups and generalised Riley slices}

\smartqed

\begin{document}

\title*{Concrete one complex dimensional moduli spaces of hyperbolic manifolds and orbifolds.}
\subtitle{{\small Kleinian groups and the Riley slice}}
\author{Alex Elzenaar \and Gaven Martin \and Jeroen Schillewaert}
\institute{Alex Elzenaar \at The University of Auckland, New Zealand, \\ \email{aelz176@aucklanduni.ac.nz}
\and Gaven Martin \at Institute for Advanced Study, Massey University, Auckland, New Zealand, \\
\email{g.j.martin@massey.ac.nz}
\and Jeroen Schillewaert \at The University of Auckland, New Zealand, \\ \email{j.schillewaert@auckland.ac.nz}
}

\maketitle

\abstract{The Riley slice is arguably the simplest example of a moduli space of Kleinian groups; it is naturally embedded in $ \CC $, and has a natural coordinate system (introduced
by Linda Keen and Caroline Series in the early 1990s) which reflects the geometry of the underlying 3-manifold deformations. The Riley slice arises in the study of arithmetic
Kleinian groups, the theory of two-bridge knots, the theory of Schottky groups, and the theory of hyperbolic 3-manifolds; because of its simplicity it provides an easy source of
examples and deep questions related to these subjects. We give an introduction for the non-expert to the Riley slice and much of the related background material, assuming only graduate
level complex analysis and topology; we review the history of and literature surrounding the Riley slice; and we announce some results of our own, extending the work of Keen and Series to the one complex dimensional moduli spaces of Kleinian groups isomorphic to $\mathbb{Z}_p*\mathbb{Z}_q$ acting on the Riemann sphere, $2\leq p,q \leq \infty$. The Riley slice is the case $p=q=\infty$ (i.e. two parabolic generators).
}


\section{Introduction}
The Riley slice is a one complex dimensional moduli space of Kleinian groups which can be embedded naturally in $ \CC $ and so is particularly amenable to study. In this article, we will define these words (Kleinian groups and moduli spaces),
outline the history of the Riley slice, and finally briefly mention some of our own results which are in the process of appearing in full.

\medskip

We have written this article in order to be accessible to the non-expert: we assume only a good understanding of areas which are plausibly familiar to most beginning graduate students (in particular, topology and complex analysis) and define almost everything else which we need from scratch. Though many of the results require only a little work to state, the study of the Riley slice has required a great deal of high-powered machinery developed in complex analysis and hyperbolic geometry in the late 20th and early 21st centuries (for instance: the theory of quasiconformal mappings and holomorphic motions, Ahlfors' null-measure conjecture, the theory of measured laminations, hyperbolic knot theory, and so forth). This is both an advantage (in that studying this rather concrete object gives an excuse to explore many of the great achievements of 20th century geometry) and a disadvantage (in that it is essentially impossible to give all detail in a short format like this article); as such we will give copious references (so that the reader inclined to explore has a variety of trails to follow) but sparse details (though we will always point towards a place where the more expert reader can look them up).

\medskip

There are two primary motivations for this work and help contextualise it. In the following loose discussion we avoid many technical details and definitions but those that matter can be found in the body of this article.  

Our first motivation is the decision problem for discrete groups, and in particular subgroups of $SL(2,\mathbb{C})$ or $\PSL(2,\mathbb{C})$.  The latter being isomorphic as a topological group to the orientation preserving  isometry group of Hyperbolic $3$-space and hence containing the fundamental groups of the most interesting class of three-dimensional geometric objects --- hyperbolic $3$-manifolds and orbifolds.  The problem is this:

\medskip 

\noindent\framebox{Decision Problem.} Given matrices $A_1,A_2,\ldots,A_n\in \PSL(2,\CC)$ decide if the group $\langle A_1,A_2,\ldots,A_n \rangle$ is discrete, and if so give a presentation for it.

\medskip

We skip a discussion of what ``decide'' means, but basically we seek a practical way to answer this question. 

If there is simply one generator,  then the answer is ``easy''.  $\langle A \rangle$ is discrete if and only if $\tr(A)\not\in (-2,2)$ or if $\tr(A)\in (-2,2)$, then $\tr(A)=2\cos(r\pi/s)$, $r,s\in \mathbb{Z}$.  If there are three or more generators,  then this problem is likely completely intractable.  This leaves the case of two generators.  The case of two generator subgroups of $\PSL(2,\mathbb{R})$ has been completely resolved by the work of Gilman \cite{gilmanAlg1,gilmanAlg2}.  Once discrete, all isomorphism types of two-generator Fuchsian groups can be found in earlier work of Purzitsky and Rosenberger, \cite{Rosenberger}.

This leaves the question of groups generated by two matrices in $\PSL(2,\CC)$. This problem goes back to to the 1950s, and we discuss it below and give some history.  But it started with the case that both generators are parabolic.  There is a good geometric reason for examining that case first as it is associated with the moduli space of the four-times punctured sphere - the simplest Riemann surface with a nontrivial Teichm\"uller space. A theorem of J\o rgensen \cite{Jorg} asserts that any subgroup of $\PSL(2,\CC)$ is discrete if and only if all of its two-generator subgroups are discrete,  a result true in much more generality, \cite{MartinActa}.  So there is utility in an effective way to determine discreteness of two-generator groups.

In a real sense the special case of two parabolic generators is now resolved. Up to conjugacy --- which preserves both discreteness and isomorphism type - there is a single complex parameter which determines the group.  The closure of the Riley slice, $ \overline{\mc{R}}$,  consists of those parameters where the group is discrete and free (there is a lot of theory behind this statement).  $ \overline{\mc{R}}$ is a closed Jordan domain (boundary a topological circle) due to Ohshika, K., Miyachi \cite{ohshika10}. A question becomes if we can determine when a complex parameter (number) lies in $\overline{\mc{R}}$. This set has a lot of additional structure, in the mid 1990s Keen and Series \cite{keen94} gave geometrically relevant coordinates for the interior of this set,  primarily through ``pleating rays". A dense set of algebraically defined curves in $\CC$ which ``land" at a dense set of the boundary $\partial \mc{R}$ (at ``cusp'' groups) by McMullen's result, \cite{mcmullen91}. Our results here give canonical open neighbourhoods of these rays which fill ${\mc{R}}$.  Using these one can construct an effective process to decide if a group is in ${\mc{R}}$ or not.  A group on the boundary $\partial  \overline{\mc{R}}$ is either a cusp group (and therefore the associated parameter is a root of a countable list of ``Farey polynomials" defined below) or is geometrically infinite with limit set the entire complex plane (and seemingly intractable computationally).  In $\CC\setminus  \overline{\mc{R}}$ the question is resolved by the recent work \cite{akiyoshi,aimi2020classification} where Aimi, Lee, Sakai, and  Sakuma, and also Akiyoshi, Sakuma, Wada, and Yamashita, resolve a conjecture of Agol.  They   identify  topologically descriptions (from which one may find a presentation)  of all quotient orbifolds for the parameters in $\CC\setminus  \overline{\mc{R}}$ for discrete groups.  These associated groups cannot be free and all these groups lie on ``extended pleating rays'' --- where a Farey word is the identity (knot and link groups) or elliptic.  Away from these parameters the groups are not discrete and their structure is quite wild \cite{martin20}.

Our work seeks to extend the programme of understanding the Riley slice to groups generated by two elements of finite order.  Such groups include all generalised triangle groups, groups with a presentation of the form 
\[   \langle a,b : a^p=b^q=R(a,b)^r =1 \rangle \]
and many groups beside.  A reasonable conjecture seems to be the following.

\begin{conjecture}\label{con2} Let $\Gamma$ be a discrete representation of $\langle a,b|a^p=b^q=1\rangle \cong \mathbb{Z}_p*\mathbb{Z}_q$ in $\PSL(2,\CC)$  and suppose $\Gamma$ is not isomorphic to $\mathbb{Z}_p*\mathbb{Z}_q$. Then there is a Farey word $w_{r/s}(a,b)\in\Gamma$ which is the identity or elliptic.
\end{conjecture}

In the case that the Farey word is the identify one would expect that the group is the orbifold fundamental group of $(p,0),(q,0)$ Dehn surgery on a two bridge link,  or, if $p=q$, a knot.  Otherwise it should follow the ``Heckoid'' structures outlined in \cite{akiyoshi,aimi2020classification}.

\medskip

This brings us to another application. Since the very beginnings of the theory of discrete groups,  mathematicians have observed deep connections between arithmetic and geometry.  Through special functions and the modular group $\PSL(2,\CC)$,  the Selberg trace formula (generally an expression for the character of the unitary representation of a Lie group $G$ on the space $L^2(G/\Gamma)$ for a discrete lattice in $G$, but initially for the traces of powers of a Laplacian defining the Selberg zeta function and the strong connections with explicit formulas in prime number theory, \cite{Selberg}.) and Margulis super-rigidity (a partial converse to the  Borel and Harish-Chandra theorem that an arithmetic subgroup in a semisimple Lie group is of finite covolume. Margulis showed that  any irreducible lattice in a semisimple Lie group of real rank larger than two is arithmetic.)

Roughly an arithmetic group is a group obtained as the integer points of an algebraic group such as $SL(n,\mathbb{Z})$.  A subgroup of $\PSL(2,\CC)$ is arithmetic if it is commensurable with another subgroup of $\PSL(2,\CC)$ isomorphic to such a group. We refer to the book of Machlachlan and Reid \cite{maclachlanreid} for a comprehensive account of the theory of arithmetic subgroups of $\PSL(2,\CC)$.  Expanding on results of Borel \cite{Borel} a necessary and sufficient criteria was determined to decide if a two generator subgroup of $\PSL(2,\CC)$ is a subgroup of an arithmetic group \cite{GMMR}. From this in \cite{maclachlan99} it was shown that there are only finitely many subgroups of $\PSL(2,\CC)$ which are either arithmetic or a subgroup of an arithmetic group which is not free,  and which is generated by two elliptic or parabolic elements. More generally we have the following conjecture.

\begin{conjecture} There are only finitely many two generator arithmetic subgroups of $\PSL(2,\CC)$.
\end{conjecture}

This might not seem too surprising as there are,  for instance,  only finitely many arithmetic subgroups of $\PSL(2,\CC)$ which lie in groups generated by reflections, \cite{agolfinite,agolfinite2}. It is not difficult to see that every two-generator group $\PSL(2,\CC)$ is a subgroup of a group generated by three involutions (elements of order two).  There seems a big difference between reflections and involutions though,  and so the conjecture remains.  

There are infinitely many arithmetic groups generated by three elements. One may see this,  for instance,  by noting the Figure of 8 knot group $K_8$ is generated by two parabolic elements, arithmetic and fibres over the circle (these sorts of things were first observed by Riley \cite{riley72,riley75}  and J\o rgensen \cite{Jorg2}) and unwrap the fibre to give a sequence of higher,  but finite, index subgroups generated by three elements and these are obviously arithmetic.

\medskip

Quite surprisingly there are only four arithmetic groups generated by two parabolic elements \cite{gehring98} and they are all knot and link complements.  A programme was instigated to identify all the arithmetic groups generated by two elliptic or parabolic elements,  for instance there are $40$ such non-uniform lattices (non-compact quotient) in $\PSL(2,\CC)$, \cite{mm3}.  This programme will soon be completed using the work described herein.  Using very recent discriminant bounds on fields with one complex place, the arithmetic criterion \cite{GMMR} identifies finitely many parameters (though sometimes many thousands).  The issue is then in deciding which of these is free on generators.  That is the primary question we give an effective way to decide.  Previous efforts in this direction tried to construct a finite volume fundamental domain for the group, \cite{mm3} or show the group in question has infinite index in a Bianchi group or other known maximal arithmetic group \cite{conder}. 

\medskip

Finally we discuss one more potential application.  In an important work, \cite{Gabai1},  Gabai gave the following very strong generalisation of the Mostow rigity theorem, in three dimensions,  provided a certain condition holds.  This is the concept of non-coalescable insulator families - structures associated to an embedded geodesic loop in a hyperbolic 3-manifold  $M$.  If $M$ has this property and if $M$ is a closed irreducible $3$-manifold and if $f:M\to N$ is a homotopy equivalence, then $f$ is homotopic to a homeomorphism.  If $f$ and $g$ are homotopic homeomorphisms of $N$, then they are isotopic and finally the space of hyperbolic metrics on $N$ is path connected.   Gabai constructs a non-coalescable insulator family for a simple closed geodesic whenever it has an embedded tubular neighbourhood in $M$ of radius larger than $\log(3)/2$.  But not every shortest geodesic has a tubular neighborhood of this radius. For example the $3$-manifold conjectured to have the third-smallest volume obtained by $(-3,2)$, $(-6,1)$-surgery on the left-handed Whitehead link.   It is believed that only  six hyperbolic 3-manifolds have shortest geodesics which do not satisfy this condition,  and for them other methods produce non-coalescable insulator families. 

The question of  $\log(3)/2$ embedded neighbourhoods is really a question about two generator groups,  with one generator of order two. To see this a geodesic will have such an embedded neighbourhood if its nearest translate (when we look at the group acting on hyperbolic $3$-space)  is further away than  $\log 3$.  Consider the group generated by the short geodesic $g$ and its nearest translate $h$.  As $g$ and $h$ are conjugate,   there is an involution $\phi$ so that $\phi g \phi =h$,  and $\langle g,\phi\rangle$ contains $\langle g,h\rangle$ with index at most two. The space of all such groups $\langle g,\phi\rangle$ up to conjugacy, discrete or not, has two complex dimensions.  Further, 
if the shortest geodesic in $M$ is quite short,  then there is always an embedded $\log(3)/2$ neighbourhood.  This is an elementary consequence of J\o rgensen's inequality. While if the shortest geodesic is long,  then its nearest translate is obviously far away too.  Thus one should explore a compact part of the moduli space of discrete groups generated by two elements,  one of which has order two,  and whose axes (the invariant hyperbolic line projecting to the geodesic representative in $M$) are relatively close.  We expect there are only finitely many such groups and very few manifold groups among them.  This is very roughly the setup for the computational exploration \cite{Gabai2} which establishes the results mentioned above in complete generality.  If there are only finitely many such exceptions in this region, then the remaining groups cannot be discrete. A contradiction to discreteness can be found by identifying a word in the group sufficiently close to the identity.  Distance to the identity is a continuous function of the parameters,  so a word sufficiently close to the identity persists, and rules out nearby groups as possibilities as well,  and therefore an open region in the search space is removed.  But how to find the words?  These words are called ``killer words'' in \cite{Gabai2} and are obtained in an ad hoc process (really a search itself).  In \cite{GMar,MM} a similar search was carried out to identify the minimal covolume lattice in $\PSL(2,\CC)$ using ``good words''.  Given that Farey words are likely relators in such discrete groups (as per Conjecture \ref{con2}) and that,  as we will describe below, the set of roots of Farey words is dense in the complement of the space of discrete free and geometrically finite groups, concrete descriptions of these root sets may well facilitate these searches.  Also clearing up the relationships between killer words, good words, Farey words and even the palindromic words (see \cite{gilmanpal}) used in discreteness criteria might yield useful insights into questions of discreteness.

\subsection{Agenda}
Broadly speaking, each section is split into two parts (which may correspond to one or more subsections): the first part of each section will be a general introduction to some
broad area which is important to the study of the Riley slice (for example, hyperbolic knot theory). Then in the second part of the section, we will explain how this area is related
to the Riley slice, giving some literature references for the details.

\begin{itemize}
  \item Section~\ref{sec:freegps} defines the notion of a Kleinian group as a discrete group of M\"obius transformations (\ref{sec:kg}) and gives the definition of the closed Riley
        slice $ \overline{\mc{R}}$ as a subsection of $ \CC $ together with the classical bounds on it (\ref{sec:rileybounds}).
  \item Section~\ref{sec:riley} describes how Kleinian groups act on hyperbolic space (\ref{sec:hyperbolicgeom}) and then gives a very fast introduction to hyperbolic knot theory (\ref{sec:knots})
        before explaining how the Riley slice appears naturally in the study of 2-bridge knots (\ref{sec:rileyhistory}).
  \item Section~\ref{sec:qc} first runs quickly through the classical theory of Riemann moduli space (\ref{sec:riemann}) before explaining the definition of a moduli space
       ---or `quasiconformal deformation space' -- of Kleinian groups motivated by the theory of holomorphic motions (\ref{sec:hm}) and then checking that the set $ \overline{\mc{R}}$
        is in fact the closure of the Riley slice as defined in Section~\ref{sec:riley} (\ref{sec:rileytop}).
  \item Section~\ref{sec:ks} explains the philosophy of the Thurston boundary of quasi-Fuchsian space and measured laminations (\ref{sec:laminations}) before briefly sketching
        the Keen--Series coordinate system of the Riley slice (\ref{sec:coordinatesystem}).
  \item Finally, Section~\ref{sec:ourwork} indicates some of our recent work which is currently appearing in detail in preprint form: we have extended the Keen--Series coordinate
        system to the so-called \emph{elliptic} Riley slices (\ref{sec:elliptic}); we have found open neighbourhoods of pleating rays in the Riley slice (\ref{sec:nbds}) and the elliptic Riley slices; and we
        have established combinatorial and dynamical results concerning pleating rays and their defining polynomials (\ref{sec:farey}).
\end{itemize}

This arrangement gives some flexibility in terms of the way in which the reader might approach this paper. Experts may wish to skim the sections introducing the Riley
slice (\ref{sec:rileybounds}, \ref{sec:rileyhistory}, \ref{sec:rileytop}, and \ref{sec:coordinatesystem}) and then read Section~\ref{sec:ourwork} in detail; students
just learning the theory of Kleinian groups may want to restrict themselves to the sections introducing general aspects of the theory and only read some of the Riley slice
material for context.

\section{Some free groups of matrices}\label{sec:freegps}
The title of this section is taken from the title of the 1955 paper of Brenner \cite{brenner55}, who (to our knowledge) was the first to explicitly
study the family of matrices which we are interested in here in the way that we are considering. We begin in Section~\ref{sec:kg} by defining Kleinian groups as Poincar\'e defined
them; then in Section~\ref{sec:rileybounds} we define the particular set of groups which we will study in the sequel.

\subsection{Kleinian groups}\label{sec:kg}
The notion of a Kleinian group was first defined by Poincar\'e in 1883 \cite{poincare83}, as a group of fractional linear transformations which acts discontinuously
on an open subset of the Riemann sphere; the name `Kleinian' was given by Poincar\'e after Klein objected to the name `Fuchsian' which Poincar\'e had previously given
to discrete subgroups of $ \Isom(\HH^2) $ \cite[p.xvi]{marden}.
We will take a slightly more general view than Poincar\'e, but before giving our definition we recall some facts from complex analysis and define
all of the terms which we have just used. Two nice introductory texts which take this classical angle are \cite{beardon} and \cite{maskit}. More modern surveys of the
field from the view of Kleinian groups (rather than, say, 3-manifolds) include those by Matsuzaki and Taniguchi \cite{matsuzakitaniguchi}, by Marden \cite{marden},
by Kapovich \cite{kapovich}, and the very nice paper by Series \cite{series05}.

We use $ \hat{\CC} $ to denote the Riemann sphere $ \CC \cup \{\infty\} $; it is well-known that if $ f : \hat{\CC} \to \hat{\CC} $ is an entire homeomorphism then
there exist $ a,b,c,d \in \CC $ satisfying $ ad - bc \neq 0 $ such that
\begin{displaymath}
  f(z) = \frac{az+b}{cz+d}
\end{displaymath}
for all $ z \in \hat{\CC} $ (with the obvious conventions $ f(-d/c) := \infty $ and $ f(\infty) := a/c $). Maps of this form are called \textit{fractional linear transformations},
or \textit{M\"obius transformations} (c.f. Chapter 3 of \cite{ahlfors}). We denote the group of all fractional linear transformations by $ \MM $.

Observe next that $ \hat{\CC} $ is identified with $ \PP\CC^1 $ (one-dimensional projective space over $ \CC $). Consider the action of the projective transformation
\begin{displaymath}
  A = \begin{bmatrix} a & b \\ c & d \end{bmatrix}\in\PSL(2,\CC)
\end{displaymath}
on $ \xi = (u,v)^t \in \PP\CC^1 $ (where we use homogenous coordinates):
\begin{displaymath}
  A\xi = \begin{bmatrix} a & b \\ c & d \end{bmatrix} \begin{bmatrix} u \\ v \end{bmatrix} = \begin{bmatrix} au + bv\\cu + dv\end{bmatrix};
\end{displaymath}
assuming $ v \neq 0 $ (so $ \xi \neq \infty $) we have that $ A $ acts on $ \CC $ as
\begin{displaymath}
  A(u/v) = \begin{bmatrix} a(u/v) + b\\c(u/v) + d\end{bmatrix} = \frac{a(u/v) + b}{c(u/v) + d}
\end{displaymath}
and so the action of $ \PSL(2,\CC) $ on $ \PP\CC^1 $ agrees with the action of $ \MM $ on $ \hat{\CC} $; just as we identify the two spheres, we identify the two groups. The matrix group $ \PSL(2,\CC) $
has an obvious topology (identify $ \SL(2,\CC) $ with a subset of $ \CC^4 $ and take the Euclidean norm; then this norm behaves well when passing to the projective quotient), and we use this topology also for $ \MM $.

\begin{definition}
  A \textit{Kleinian group} is a discrete subgroup of $ \MM $.
\end{definition}

\begin{figure}
  \centering
  \includegraphics[width=0.7\textwidth]{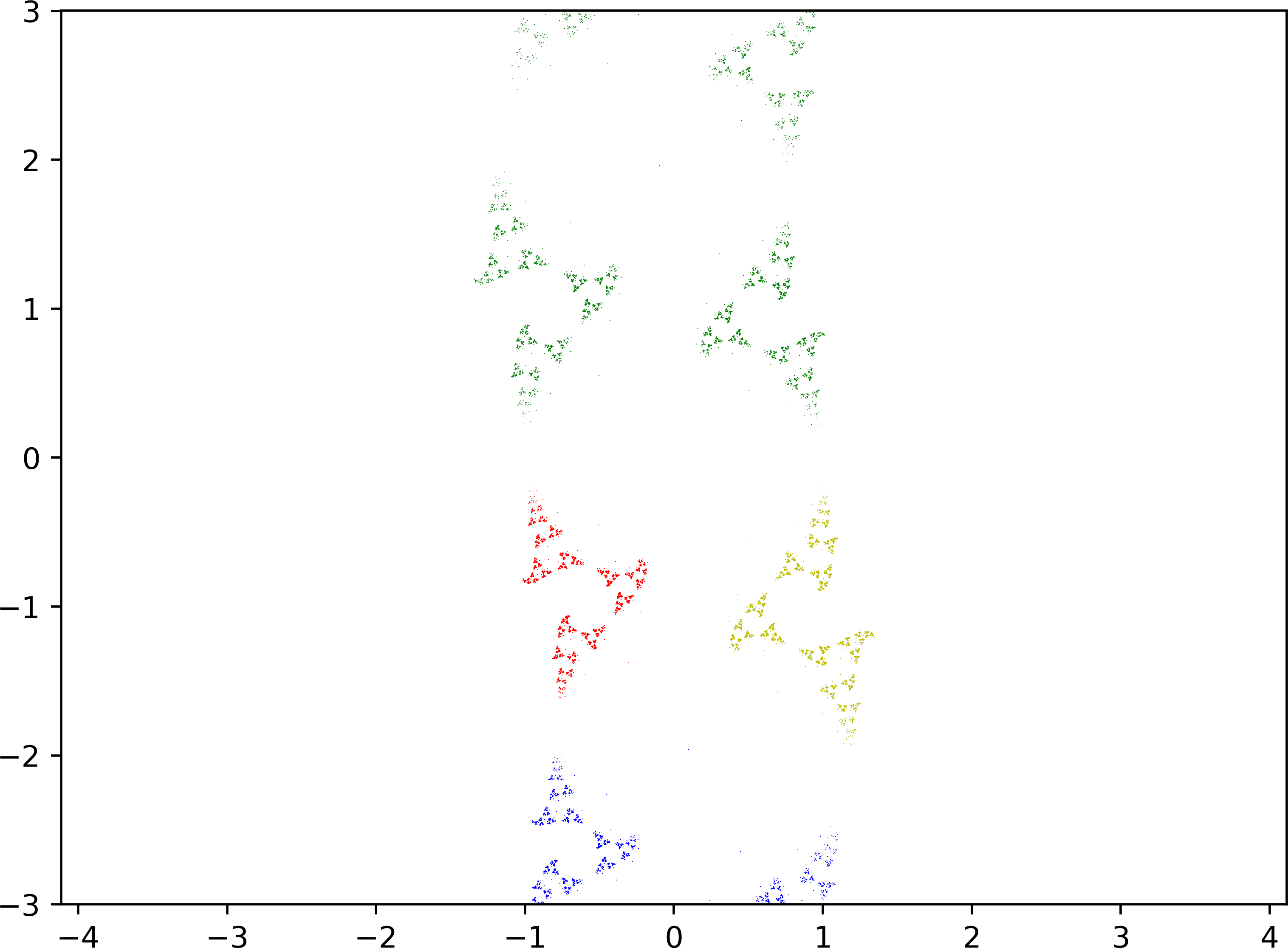}
  \caption{The limit set of a Kleinian group, generated by two matrices with respective traces $ 2\pm 0.2i $ (c.f. Figure 8.5(iv) of \cite{indras_pearls}).\label{fig:indra}}
\end{figure}

We now define the term `discontinuous' which we used in giving Poincar\'e's definition of a Kleinian group, and explain why our definition is slightly more general. Let $ G $ be a
group acting by homeomorphisms on a topological space $ X $. We can view this action as a dynamical object, and ask for the orbits of points under the group
action. We say that $ x \in X $ is a \textit{limit point} of $ G $ if there is some $ u \in X $ such that $ x $ is an accumulation point of the orbit $ Gu $;
the set of all limit points is denoted by $ \Lambda(G,X) $. In our case, the set $ X $ will always be $ \hat{\CC} $ and so we simply write $ \Lambda(G) $
for the limit set of a group of M\"obius transformations $ G $. Directly from its definition the limit set is always a closed set. We exhibit an example as Fig.~\ref{fig:indra}.

For a proof of the following proposition, see Section~II.E of \cite{maskit}.
\begin{proposition}
  Suppose that $ G \leq \MM $. The following two conditions are equivalent for $ x \in \hat{\CC} $:
  \begin{enumerate}
    \item $ x \in \hat{\CC} \setminus \Lambda(G) $
    \item There exists a neighbourhood $ U $ of $ x $ such that $ U \cap gU \neq \emptyset $ for only finitely many $ g \in G $. \qed
  \end{enumerate}
\end{proposition}
If $ x $ satisfies condition (2), then we say that $ G $ acts \textit{discontinuously} at $ x $; the set of all points at which $ G $ acts
discontinuously, which is identical to the open set $ \hat{\CC}\setminus\Lambda(G) $ by the proposition, is denoted by $ \Omega(G) $ and is called
the \textit{ordinary set} of $ G $. If $ G $ is a Kleinian group and $ \Omega(G) $ is non-empty, then $ \Omega(G)/G $ naturally has a complex structure;
if $ G $ is finitely generated, by Alhfors' finiteness theorem the space $ \Omega(G)/G $ is a finite union of Riemann surfaces, each with finitely many
punctures and no other boundary components \cite{ahlfors64,ahlfors64err} (for a modern proof and references to several other proofs see \cite[\S 8.14]{kapovich}).

Poincar\'e's definition of a Kleinian group was a group $ G $ such that $ \Omega(G) \neq \emptyset $;
it is not too hard to prove that if $ G \leq \MM $ acts discontinuously at a point $ x \in \hat{\CC} $, then $ G $ is discrete (one uses
the fact that a subgroup of $ \MM $ is discrete iff there are no nontrivial sequences in the group converging to the identity, see \cite[II.C.3]{maskit});
on the other hand, there are discrete groups of M\"obius transformations which act discontinuously nowhere \cite[VIII.G]{maskit}, so our definition is strictly more general
than Poincar\'e's. (The advantage of our definition is that it includes groups which act discontinuously on the half-space model of hyperbolic space whose sphere at infinity is
$ \hat{\CC} $, even if they do not act discontinuously at infinity; we will take this story up again in Section~\ref{sec:hyperbolicgeom}.)

As we outlined in the introduction,  in general, the question of deciding whether an arbitrary subgroup $ G \leq \MM $ is discrete (as well as the related problems of computing $ \Omega(G) $ and $ \Lambda(G) $ explicitly) is a hard one. In the next
section, we will describe the particular space of groups of M\"obius transformations which we are interested in; we will see that the boundary between the set of discrete groups and the set of non-discrete groups
in this space is highly intricate, even though the space is arguably the simplest possible space of groups for which we might ask this question.

\subsection{The space of discrete, free groups generated by  $2$-parabolics}\label{sec:rileybounds}
\begin{figure}
  \centering
  \includegraphics[width=0.7\textwidth]{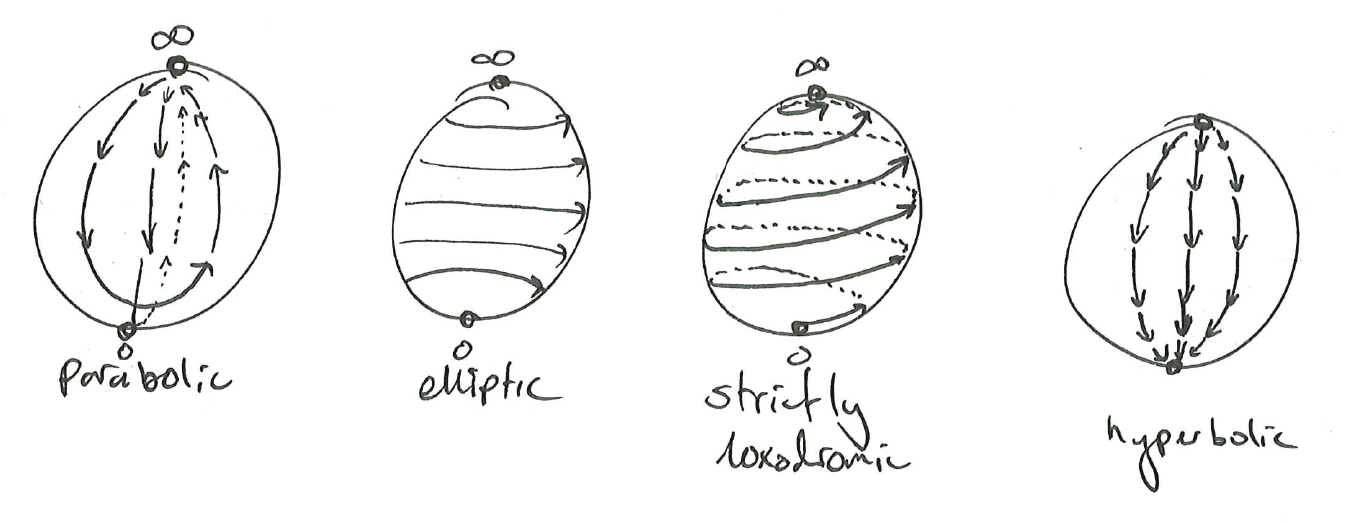}
  \caption{The flow of points on $ \hat\CC $ under the four kinds of elements of $ \MM $.\label{fig:conjugacy_classes}}
\end{figure}

We return briefly to the subject of orbits of a group of M\"obius transformations.
\begin{proposition}
  For $ \lambda \in \CC\setminus\{0\} $, define the function $ f_{\lambda^2} $ by $ f_{\lambda^2}(z) := \lambda^2 z $. Every M\"obius transformation is
  conjugate by an element of $\MM$ either to $ z \mapsto z+1 $, or to exactly one of the functions $ f_{\lambda^2}$. \qed
\end{proposition}
This gives a taxonomy of M\"obious transformations. Let $ g \in \MM $ be non-trivial,  that is $g$ is not the identify; then $ g $ is called
\begin{itemize}
  \item \textit{parabolic}, if it is conjugate to the translation $ z \mapsto z + 1 $ (equivalently, $ \tr^2 g = 4 $);
  \item \textit{elliptic}, if it is conjugate to some rotation $ z \mapsto \lambda^2 z $ with $ |\lambda| = 1 $ (equivalently, $ \tr^2 g \in [0,4) $);
  \item \textit{hyperbolic}, if it is conjugate to some dilatation $ z \mapsto \lambda^2 z $ with $ \lambda \in \RR $ and $ |\lambda| \neq 1 $ (equivalently, $ \tr^2 g \in (4,\infty) $);
  \item \textit{strictly loxodromic}, if it is conjugate to $ z \mapsto \lambda^2 z $ with $ \lambda \in \CC \setminus \RR $ and $ |\lambda| \neq 1 $ (equivalently, $ \tr^2 g \in \CC\setminus[0,\infty) $).
\end{itemize}
Parabolic elements have exactly one fixed point on the Riemann sphere $ \hat\CC $; all other non-identity elements have exactly two fixed points, and in the case of hyperbolic
and strictly loxodromic elements (which together form the class of \textit{loxodromic} elements) these fixed points form a source-sink pair (Fig.~\ref{fig:conjugacy_classes}).

In 1947, Sanov \cite{sanov47} gave the following example of a free Kleinian group:
\begin{displaymath}
  \left\langle \begin{bmatrix} 1 & 2 \\ 0 & 1 \end{bmatrix}, \begin{bmatrix} 1 & 0 \\ 2 & 1 \end{bmatrix} \right\rangle.
\end{displaymath}
This group was of interest because it was a free subgroup of $ \PSL(2,\ZZ) $; previous work had found free subgroups of $ \MM $, but all
of these had transcendental entries. In 1955, Brenner \cite{brenner55} extended this result, showing that whenever $ |m| \geq 2 $,
\begin{displaymath}
  G_m = \left\langle A_m = \begin{bmatrix} 1 & m \\ 0 & 1 \end{bmatrix}, B_m = \begin{bmatrix} 1 & 0 \\ m & 1 \end{bmatrix} \right\rangle.
\end{displaymath}
is free. Observe that both generators of this group are parabolic (since they both have trace 2). Let $ T \in \PSL(2,\CC) $ be a matrix conjugating $ A_m $ to the matrix
\begin{displaymath}
  \begin{bmatrix} 1 & 1 \\ 0 & 1 \end{bmatrix}
\end{displaymath}
which represents $ z \mapsto z + 1 $; for example, we can choose
\begin{displaymath}
  T = \begin{bmatrix} \sqrt{m} & 0 \\ 0 & \frac{1}{\sqrt{m}} \end{bmatrix}
\end{displaymath}
and then compute that
\begin{displaymath}
  T^{-1}G_mT = \left\langle \begin{bmatrix} 1 & 1 \\ 0 & 1 \end{bmatrix}, \begin{bmatrix} 1 & 0 \\ m^2 & 1 \end{bmatrix} \right\rangle.
\end{displaymath}

There is a famous inequality due to Shimizu \cite{shimizu63} and Leutbecher \cite{leutbecher67} (see also \cite[II.C.5]{maskit}) which gives a necessary condition for a group of this form to be discrete:
\begin{lemma}[Shimizu--Leutbecher lemma]
  Suppose $ A,B \in \PSL(2,\CC) $ are of the form
  \begin{displaymath}
    A = \begin{bmatrix} 1&1\\0&1\end{bmatrix}, B = \begin{bmatrix} a & b\\c & d \end{bmatrix}.
  \end{displaymath}
  If $ \langle A,B \rangle $ is discrete, then either $ c = 0 $ or $ |c| \geq 1 $. \qed
\end{lemma}

Combining this with the result of Brenner, we obtain bounds on the set of $ \rho \in \CC $ such that
\begin{displaymath}
  \Gamma_\rho := \left\langle X = \begin{bmatrix} 1 & 1 \\ 0 & 1\end{bmatrix}, Y_\rho = \begin{bmatrix} 1 & 0 \\ \rho & 1 \end{bmatrix} \right\rangle
\end{displaymath}
is free and discrete. Let $ \overline{\mc{R}} $ be this set:
\begin{lemma}\label{lem:crudebound}
  Let $ \rho \in \CC $. If $ \rho \in \overline{\mc{R}} $, then $ |\rho| \geq 1 $. On the other hand, if $ |\rho| \geq 4 $, then $ \rho \in \overline{\mc{R}} $. \qed
\end{lemma}

\begin{figure}
  \centering
  \includegraphics[width=0.5\textwidth]{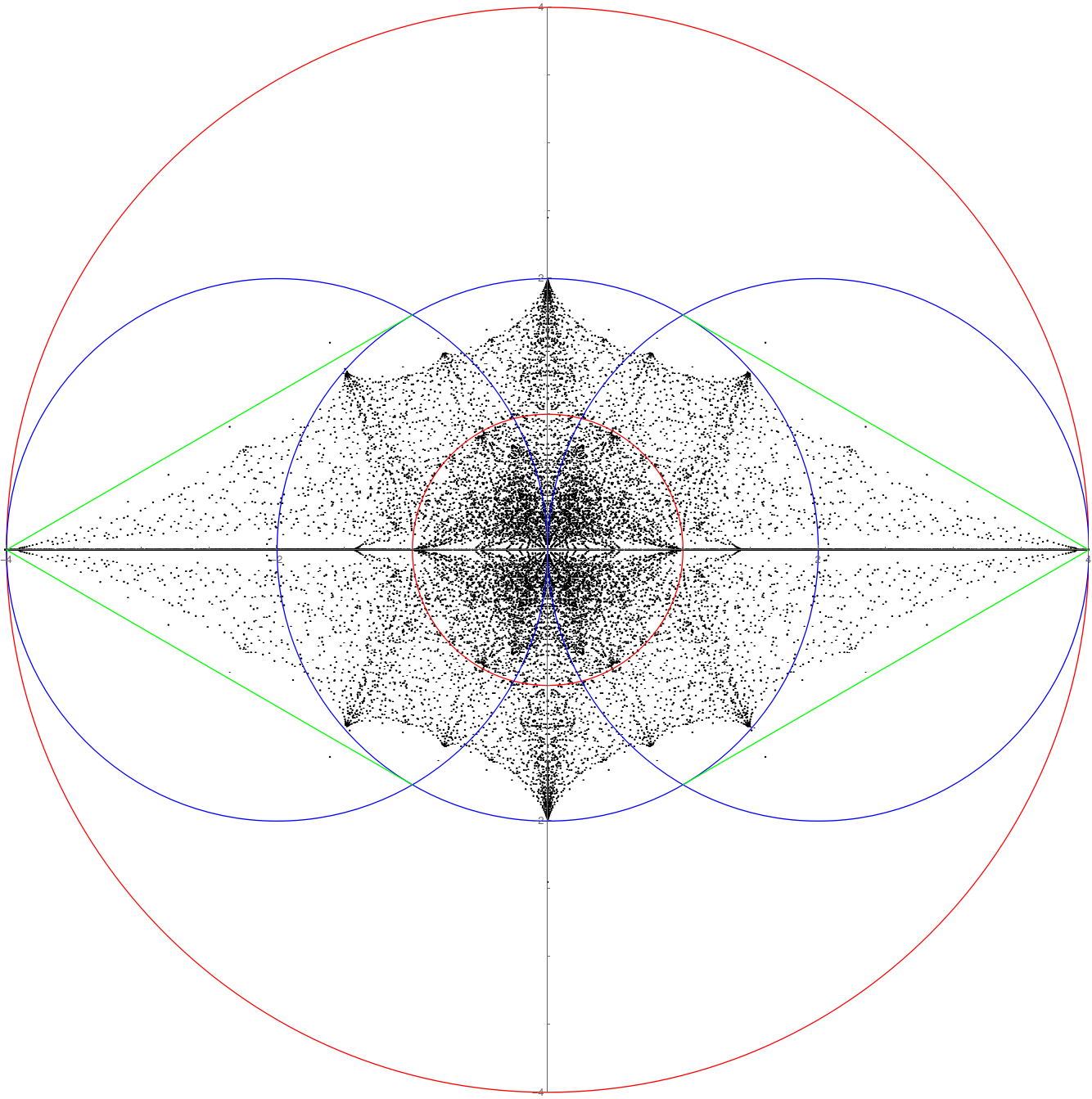}
  \caption{Crude bounds on the Riley slice.\label{fig:crudebound}}
\end{figure}

We will call the set $ \overline{\mc{R}} $ the \textit{closed Riley slice} (we will explain this terminology in the next section). Using techniques which we will reference
in Section~\ref{sec:ks} and Section~\ref{sec:ourwork}, we are able to draw a picture of this set which shows that it has a very intricate structure which is not captured by
the crude bounds of Lemma~\ref{lem:crudebound}; the closed Riley slice is the region of $ \CC $ which is outside the shaded `eye' of Fig.~\ref{fig:crudebound}, and the bounds of
Lemma~\ref{lem:crudebound} are the two red circles.

In 1958, Chang, Jenning, and Ree \cite{chang58} improved the bound on the interior of the Riley slice from `$ |\rho| \geq 4 $' to `$\rho$ lies outside the three circles
of radius $2$ centred at $ -2 $, $ 0 $, and $ 2 $ respectively'; this is shown in blue in Fig.~\ref{fig:crudebound}. A couple of years later, Ree \cite{ree61} showed
various density results for the complement of the Riley slice: for instance, that non-free groups are dense in the closed unit disc $ \{\rho\in\CC:|\rho|\leq1\}$;
and that the two open intervals $ (-4,4) $ and $ (-2i,2i) $ are contained in open sets in which non-free groups are dense.

Finally, in 1969, Lyndon and Ullman \cite{lyndon69} improved the interior bound once more; this bound is shown in green on Fig.~\ref{fig:crudebound}.
\begin{theorem}[Lyndon and Ullman]
  Let $ LU $ be the convex hull in $ \CC $ of the circle $ \{z\in\CC:|z|=2\} $ and the two points $ \pm 4 $. Then $ \overline{\mc{R}} $ is contained
  in $ \overline{\CC\setminus LU} $. \qed
\end{theorem}
For some further study of the structure of the Riley slice boundary using Lyndon and Ullman's result, see Section~4.2 of \cite{ems21}. Beardon has given improvements of
some of the other results of Lyndon and Ullman on the set of free groups of the form $ \Gamma_\rho $ \cite{beardon93}.

Some more modern bounds on sets intimately related to the Riley slice as defined here (and the elliptic Riley slice which we will describe in Section~\ref{sec:ourwork}) are
given in \cite{martin20b}. The non-discrete groups of the form $ \Gamma_\rho $ have also been studied, for example by Beardon \cite{beardon96} and Martin \cite{martin20};
the non-free groups have been studied by Bamberg \cite{bamberg00}; and the arithmetic groups (whose definition would take us too far afield, see the elementary book by Voight \cite{voight}
and the more difficult book by Maclachlan and Reid \cite{maclachlanreid}) have been studied by Gehring, Maclachlan, and Martin \cite{gehring98}, see also the discussion and citations in
the introduction to \cite{ems21}. Work on parabolic conjugacy classes in the groups $ \Gamma_\rho $, whether free or not, was done by Gilman \cite{gilman08} in an attempt to
explain the patterns of density seen in pictures of the Riley slice exterior like those of Fig.~\ref{fig:crudebound} above (which was drawn, following Keen and Series \cite{keen94}, by
drawing the set of groups where certain \textit{Farey words} become accidentally parabolic---we discuss this more in Section~\ref{sec:ks} below).

\section{Riley's study of knot groups}\label{sec:riley}
In this section we study hyperbolic geometry, and show that the Riley slice appears naturally in association with representations of the fundamental groups
of 2-bridge knot groups.

\subsection{Hyperbolic geometry}\label{sec:hyperbolicgeom}
We do not give an account here of the interesting history of `non-Euclidean geometry'; some reading on this matter
includes the first chapter of \cite{coxeterNE} as well as the monograph \cite{rosenfeld}. We simply note that hyperbolic
geometry was first studied as the geometry in which Euclid's \textit{Postulate V}
\begin{quote}
  For every point $ A $, and any line $ r $ not through $ A $, there is not more than one line through $ A $, in the plane spanned by $ A $ and $ r $, which does not meet $ r $
\end{quote}
\noindent (here formulated as in \cite[\S16.1]{coxeterI}) is replaced by the postulate
\begin{quote}
  For every point $ A $, and any line $ r $ not through $ A $, there is more than one line through $ A $, in the plane spanned by $ A $ and $ r $, which does not meet $ r $,
\end{quote}
\noindent and in this form the first people to study it seriously were Lobačevskiĭ \cite{lobacevskii29} and J. Bolyai \cite{bolyai32}.

\begin{figure}
  \centering
  \includegraphics[width=0.5\textwidth]{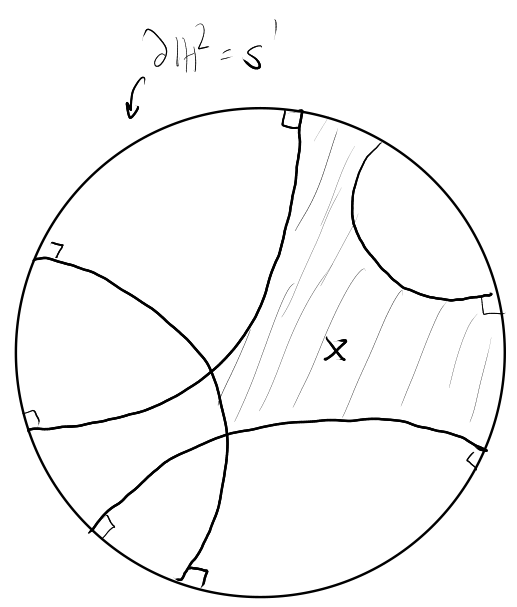}
  \caption{A polygon $ X $ in the ball model of $ \HH^2 $, with four edges given by hyperbolic geodesic segments and two edges on the sphere $ S^1 $ at infinity.\label{fig:hyper2}}
\end{figure}

For us, \textit{hyperbolic $ n$-space} is the upper half-plane $ \HH^n := \{ (x_1,...,x_n) \in \RR^n : x_n > 0 \} $
endowed with the metric given by the differential
\begin{displaymath}
  ds^2 = \frac{dx_1^2 + \cdots + dx_n^2}{x_n^2};
\end{displaymath}
geodesics in this metric are the intersections with $ \HH^n $ of lines and circles in $ \RR^n $ orthogonal to the plane $ x_n = 0 $. This space is isometric (after taking a sphere inversion through a sphere tangent
to $ x_n = 0 $, translating, and scaling) to the space consisting of the interior of the $ n$-disc $ \BB^n := \{ x \in \RR^n : |x|<1 \} $ with the metric
\begin{displaymath}
  ds^2 = \frac{2dx_1^2 + \cdots + 2dx_n^2}{(1-dx_1^2-\cdots-dx_n^2)^2},
\end{displaymath}
in which geodesics are circle arcs orthogonal to the boundary $ S^{n-1} $ (see Fig.~\ref{fig:hyper2}).

We call $ \HH^n $ the \textit{half-space model} and $ \BB^n $ the \textit{ball model} of hyperbolic $ n$-space; these are two views of the unique Riemann manifold of constant
curvature $ -1 $. Some introductory texts on hyperbolic geometry include \cite[Chapter 16]{coxeterI} (which takes an axiomatic view), \cite[Chapter 7]{beardon} (which takes an analytic
view like that we take here), and \cite[Chapters 1--5]{ratcliffe} (which starts with the point of view of Lorentzian space and proceeds to discuss both axiomatic and geometric views).

Clearly both models of hyperbolic space have some kind of `boundary': this is a common feature of metric spaces of non-positive curvature \cite[Chapters 8 and 9]{bridson_haefliger} and in
this case, the boundary of hyperbolic space is a sphere (either $ \RR^{n-1} \cup \{\infty\} $ for the half-plane model, or the usual unit sphere $ S^{n-1} $ for the ball model). We will
primarily be interested in hyperbolic 3-space, in which case we can naturally identify $ \RR^2 $ with $ \CC $ and hence the sphere at infinity with the Riemann sphere $ \hat{\CC} $.

Let $ \phi $ be an (orientation-preserving) isometry of $ \HH^3 $. Then $ \phi $ extends naturally to a conformal homeomorphism (hence M\"obius transformation) on the boundary $\partial   \HH^3  =\hat{\CC}$---this is a non-trivial
observation but isn't too hard to check, see \cite[IV.A]{maskit}---and in fact this map $ \Isom(\HH^3) \to \MM $ is invertible (which is not a feature of all metric spaces of non-positive curvature). The inverse
map is the so-called \textit{Poincar\'e extension}, where to define the action of a conformal map $ \psi $ on the boundary $ \hat\CC $ on a point $ x \in \BB^n $ one first finds a pair of
intersecting geodesics through $ x $, computes the images of the endpoints at infinity of these geodesics under $ \psi $, and uses these to draw two new geodesics which will intersect at
the image of $ x $ under the extension of $ \psi $ to $ \HH^3 $ (see the final paragraph of \cite[A.13]{maskit}, or for more detail \cite[\S3.3]{beardon} and \cite[\S4.4]{ratcliffe}). In any case, once we accept
that there is a natural geometric bijection between $ \MM $ and the group of isometries of $ \HH^3 $ (which in fact is compatible with the compact-open topology on the isometry group) then we may
give an alternative definition of a Kleinian group: a discrete group of isometries of $ \HH^3 $.

We now follow Thurston \cite{thurstonB} (see also \cite[Chapter 8]{ratcliffe} and \cite{benedetti_petronio}) to define \textit{hyperbolic 3-manifolds}. A hyperbolic 3-manifold is a Hausdorff second-countable
topological space $ M $ equipped with a \textit{hyperbolic atlas}: an open cover $ U_{\alpha} $ with topological embeddings $ \phi_\alpha : U_\alpha \to \HH^3 $ ($ \alpha $ in some index set $ A $)
where every transition map $ \phi_\beta \phi_\alpha^{-1} $ ($\alpha,\beta\in A$) is the restriction of some isometry of $ \HH^3 $ to the domain of definition.
The holonomy group of $ M $ is a structure which detects the tiling that is induced when $ M $ is rolled out onto its universal cover $ \HH^3 $ in a way that respects the geometry (this process is
called \textit{developing}; see \cite[\S 3.4]{thurstonB} for an intuitive description, much more detail may be found in \cite[\S 8.4]{ratcliffe}). If $ \Gamma $ is the holonomy group
of a connected and metrically complete hyperbolic 3-manifold $ M $, then it is a Kleinian group and $ M $ is isometric to $ \HH^3/\Gamma $ \cite[Theorem 8.5.9]{ratcliffe}. Conversely, if $ \Gamma $ is
a Kleinian group then $ \HH^3/\Gamma $ is a hyperbolic manifold if $ \Gamma $ is torsion-free; if $ \Gamma $ is allowed to have torsion then the quotient is a hyperbolic orbifold (i.e. a hyperbolic manifold
with possible cone points---see Chapter~13 of Thurston's notes \cite{thurstonN}, Chapter~13 of \cite{ratcliffe}, Chapter~6 of \cite{kapovich}, and Chapter~III.$\mc{G}$ of \cite{bridson_haefliger}).

Hyperbolic 3-manifolds are a special class of so-called \textit{geometric manifolds}. These are manifolds with atlases that give a structure locally modelled on a metric space $ X $ which satisfies
five axioms generalising Euclid's five (written below, verbatim from \cite[Section 1.2]{coxeterI}, in italics). This generalisation is described in much more detail in Chapter 8 of \cite{ratcliffe},
and a more differentio-geometric definition of geometric metric spaces is given by Thurston in Section 3.8 of \cite{thurstonB}.

\begin{enumerate}
  \item \emph{A straight line may be drawn from any point to any other point.} For any pair $ x,y \in X $ there exists a geodesic segment
        joining $ x $ to $ y $.
  \item \emph{A finite straight line may be extended continuously in a straight line.} Given a geodesic segment $ \gamma : [0,1] \to X $
        there exists a unique geodesic line $ \hat\gamma : \RR \to X $ with $ \hat\gamma\restriction_{[0,1]}=\gamma $.
  \item \emph{A circle may be described with any centre and any radius.} There exists a continuous function $ \sigma : \RR^n \to X $ and a real
        number $ \varepsilon > 0 $ such that $ \sigma $ is a homeomorphism of $ B(0,\varepsilon) $ onto $ B(\sigma(0),\varepsilon) $, such that
        for any $ x \in S^{n-1} $ the map $ \gamma : \RR \to X $ defined by $ \gamma(t) = \sigma(tx) $ is a geodesic line such that $ \gamma\restriction_{[-\varepsilon,\varepsilon]} $
        is a geodesic segment.
  \item \emph{All right angles are equal to each other.} $ X $ is homogeneous (i.e. for all $ x, y \in X $ there exist isometric neighbourhoods $ U $ of $ x $ and $ V $ of $ y $).
\end{enumerate}

There are eight such geometries in three dimensions: Euclidean space $ \mathbb{E}^3 $, spherical space $ S^3 $, and hyperbolic space $ \HH^3 $; the two product
geometries $ S^2 \times \mathbb{E}^1 $ and $ \HH^2 \times \mathbb{E}^1 $; the universal cover of $ \SL(2,\RR) $; solvegeometry $ \mathsf{Sol} $; and nilgeometry $ \mathsf{Nil} $.
Some excellent visualisation work has been done with these geometries; for instance, the work of Coulon, Matsumoto, Segerman, and Trettel \cite{cmstGH,cmstPP}

It is a great insight of Thurston that `most' geometric manifolds are hyperbolic. He established this in many interesting cases, particularly the Haken case (compact, orientable, irreducible $3$-manifold that contains an orientable, incompressible surface) and promoted the ``Geometrisation Conjecture'' which was completely solved by Grigori Perelman,  see the book  \cite{Morgan},  using ideas of Richard S. Hamilton who developed the the Ricci flow.    We recall that if $S$ is  a properly embedded surface in a $3$-manifold $M$,  then a compression disc $D$ (for $S$) is an embedded disk in $M$ with $D \cap S =\partial D$, but the loop $\partial D\subset S$ does not  bound a disc in $S$. If no such compression disk exists, then $S$ is is said to be incompressible.

\subsection{Knots}\label{sec:knots}
The theory of knots was found to be related to hyperbolic geometry by Robert Riley and William Thurston in the mid 1970s. In this section, we define
some of the basic concepts from geometric knot theory and explain exactly how hyperbolic geometry comes into it all. The three most useful references
for the parts of knot theory which are of interest to our story are (in rough order of expected prior knowledge) the books by Purcell \cite{purcell}, Rolfsen \cite{rolfsen},
and Burde and Zieschang \cite{burde}. Unfortunately there seems not to be a detailed monograph on the history of knot theory; \cite{burde} has some historical notes at the end of each chapter,
and the book chapter by Epple \cite{histoftop}, though focused primarily on the history of knot invariants, does have some discussion of the history of knot geometry. In addition, Chapter~3 of
\textit{In the tradition of Thurston} \cite{thurstonTrad} has a detailed description of Thurston's contributions as well as some earlier history.

\begin{definition}
  A \textit{link} is a topological embedding of a collection of disjoint circles into $ S^3 $; the image of each circle is called a \textit{component}
  of the link, and a \textit{knot} is a link with a single component. We restrict ourselves to links which are \textit{tame}: every component must be ambient
  isotopic in $ \RR^3 $ to a simple closed piecewise linear curve. In general, we view two links as equivalent if they are ambient isotopic in $ S^3 $.
\end{definition}

Clearly the topological type of the complement of a knot $ k $ in $ S^3 $ (formally, $ \overline{S^3 \setminus V} $ where $ V $ is an open tubular neighbourhood of $k$) is an isotopy invariant of the knot.
That the converse is also true is a famous result of Gorden and Luecke \cite{gordon89}. It follows that studying the topology of knot complements is equivalent to studying the structure of knots. (The corresponding
result is not true for links.) It is an insight of Thurston's programme of geometry that one can study the topology of spaces by studying their geometry.

In 1982, Thurston announced the following pair of results, Theorems~\ref{thm:thurston1} and~\ref{thm:thurston2} in his seminal paper \cite{thurston82}.
\begin{theorem}[{Corollary 2.5 of \cite{thurston82}}]\label{thm:thurston1}
  If $ k \subseteq S^3 $ is a knot, $ S^3 \setminus k $ has a geometric structure (i.e. admits a geometric atlas) iff $ k $
  is not a satellite knot. It has a hyperbolic structure iff, in addition, $k$ is not a torus knot.\qed
\end{theorem}
Thurston conjectured Theorem~\ref{thm:thurston1}, which essentially shows that most knot complements admit a hyperbolic structure, based on the work of Riley which we discuss in the next section.
The history surrounding this discovery is very interesting; various accounts beyond \cite{thurston82} include \cite{Riley13} and the accompanying commentary \cite{Brin13}, and the
additional references given in the historical notes to Section 10.3 on p.504 of \cite{ratcliffe}.

If $ L $ is a link in a 3-manifold $ M $, then a \textit{Dehn surgery} along $ L $ is effected by cutting out a tubular neighborhood of $ L $, performing a diffeomorphism on each torus component of that
neighbourhood, and gluing it back into $ M $ (see Section~1.5 of \cite{ratcliffe}, Chapter~6 of \cite{purcell}, Section E.4 of \cite{benedetti_petronio}).
\begin{theorem}[{Theorem 2.6 of \cite{thurston82}}]\label{thm:thurston2}
  Suppose $ L $ is a link in a 3-manifold $ M $ such that $M\setminus L$ has a hyperbolic structure.
  Then most manifolds obtained from $M$ by Dehn surgery along $L$ have hyperbolic structures.
\end{theorem}

Every closed 3-manifold is obtained from $ S^3 $ by Dehn surgery along some link whose complement is hyperbolic; this follows from the Lickorish-Wallace theorem
(which shows that every manifold comes from Dehn surgery along some link) \cite{wallace60,lickorish62} (see also \cite[\S 9I]{rolfsen} and Theorem~6.5 of \cite{purcell}) together with
the observation that every link in the 3-sphere is a sublink of a hyperbolic link \cite{myers82}: just do the trivial Dehn surgery on the added components. Thus Theorem~\ref{thm:thurston2}
shows that most closed 3-manifolds admit hyperbolic structures.

\subsection{Riley's study of 2-bridge knots and links}\label{sec:rileyhistory}
The Riley slice arises naturally in the study of knots. In particular, a dense subset of the groups in the Riley slice can be viewed as degenerations of holonomy groups of 2-bridge knot
complements. In this section we briefly discuss Riley's study of these groups, and explain how to obtain the Riley slice from them. Much more detail may be found in Section~2 of our upcoming
preprint, \cite{ems22}, as well as in the papers of Akiyoshi, Ohshika, Parker, Sakuma, and Yoshida \cite{akiyoshi2020classification} and Aimi, Lee, Sakai, and Sakuma \cite{aimi2020classification}.

We first define 2-bridge links, following the treatment of \cite[Chapter 10]{purcell}. A \textit{rational tangle} is obtained by drawing two disjoint arcs on the boundary of a 3-ball $ B $ and then
pushing the interiors of the arcs slightly into the 3-ball, forming a pair of disjoint arcs embedded in $ B $ with endpoints on $ \partial B $ so the boundary is a sphere with four marked points.
A \textit{2-bridge link} is a link obtained by joining these marked points up with a pair of disjoint arcs on the boundary (the \textit{bridge arcs}). There is a natural bijection due to
Schubert \cite{schubert56}, described in detail as Algorithm~2.1 of \cite{ems22}, which sets up a bijection between the set of 2-bridge links (up to ambient isotopy) and $ \hat{\QQ} := \QQ \cup \{\infty\} $;
the rational number associated to a link is called the \textit{slope} of the link.\footnote{Our `slope' is actually the \emph{reciprocal} of Schubert's normal form.}

One of the simplest examples of a 2-bridge link is the figure 8 knot $ k $ (which has slope $ 3/5 $). It is important historically as it was the first known example of a link
whose complement manifold admits a hyperbolic structure; it is Riley's proof of this which (according to Thurston \cite{thurston82}) inspired the conjecture
of Theorem~\ref{thm:thurston1}.

Riley's proof of the hyperbolicity is described in \cite{riley75b}; it is a consequence of his algebraic theory of parabolic representations of 2-bridge link groups developed in \cite{riley72,riley75}.
The main idea is to consider a representation $ \theta : G \to \PSL(2,\CC) $ (where $ G = \langle x, y : x^{-1}yxy^{-1} x = yx^{-1}yxy^{-1} \rangle $ is the fundamental group of the knot complement) defined by
\begin{displaymath}
  \theta x := A = \begin{bmatrix} 1 & 1  \\ 0 & 1 \end{bmatrix} \qquad \theta y := B = \begin{bmatrix} 1 & 0 \\ -\omega & 1 \end{bmatrix}
\end{displaymath}
where $ \omega $ is a primitive cube root of unity. Then:
\begin{theorem}[Theorem 1 of \cite{riley75b}]\label{thm:riley}
  If $ \Gamma := \langle \theta x, \theta y \rangle $ then $ \Gamma $ has presentation
  \begin{equation}\label{eqn:overpres}
    \langle A, B : WAW^{-1} = B \rangle
  \end{equation}
  where $ W = A^{-1} BAB^{-1} $, and so $ \theta $ gives an isomorphism $ G \simeq \Gamma $. \qed
\end{theorem}
Hyperbolicity of the complement manifold (Corollary on p.284 of \cite{riley75b}) now follows from this theorem and some standard results on 3-manifold topology.

Every 2-bridge link has a fundamental group with an presentation on two generators like (\ref{eqn:overpres}) for the figure 8 knot (such a representation is called
an \textit{over-presentation}---the two generators correspond to loops about the two bridge arcs), and Riley showed that every such group admits a non-abelian representation
with the images of the two generators being parabolic (Theorem 2 of \cite{riley72}). In fact, in 2002 Agol \cite{agol02} announced a full classification of the non-free discrete groups generated by
two non-parabolic elements (later given a formal and correct presentation in \cite{akiyoshi2020classification} and \cite{aimi2020classification}, see also Section~2.2 of \cite{ems22}).

\begin{figure}
  \centering
  \includegraphics[width=\textwidth]{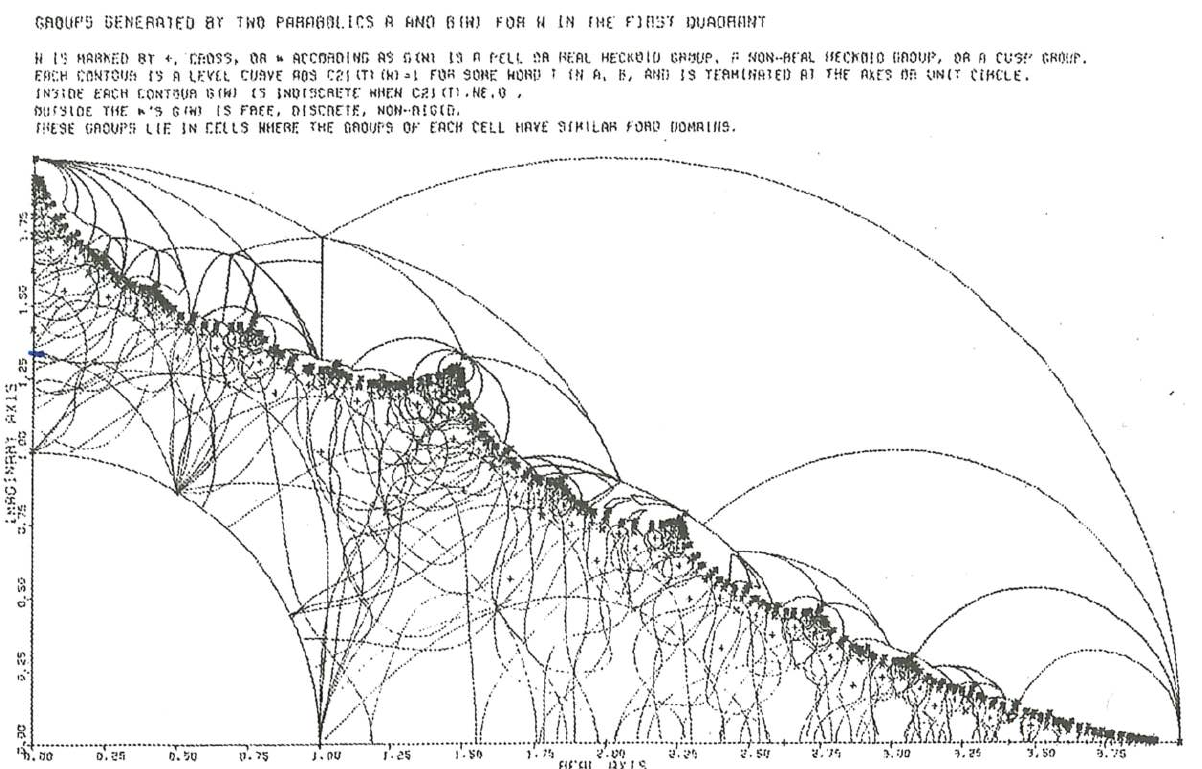}
  \caption[Riley's plot of the Riley slice]{Riley's plot of 2-bridge link groups in the $(+,+)$-quadrant of $ \CC $.\label{fig:riley_orig_pic}}
\end{figure}

If $ X $ and $ Y $ are two parabolics which do not share fixed points then, up to conjugacy, they are represented by matrices
\begin{displaymath}
  X = \begin{bmatrix} 1 & 1 \\ 0 & 1 \end{bmatrix}, Y_\rho = \begin{bmatrix} 1 & 0 \\ \rho & 1 \end{bmatrix}
\end{displaymath}
where $ \rho \in \CC $. This normalisation associates a complex number $ \rho $ with each 2-bridge link group representation. Riley devoted much effort to finding the set of $ \rho \in \CC $ such that
\begin{displaymath}
  \Gamma_\rho := \langle X, Y_\rho \rangle
\end{displaymath}
is a 2-bridge link group, producing the famous image in Figure~\ref{fig:riley_orig_pic}.

The most obvious feature of Figure~\ref{fig:riley_orig_pic} is the fractal curve running diagonally across the centre of the frame. This is the set of so-called \textit{cusp groups}.
The half of the plot on the left (closer to 0) contains (marked with $ + $ or $ \times $) various \textit{Heckoid groups} (introduced by Riley in \cite{riley92} and formalised by Lee
and Sakuma in \cite{lee13}), which are the non-free Kleinian groups mentioned above in association to Agol's 2002
announcement. Observe that these points appear to lie on curves (called \textit{extended pleating rays}) which end on the diagonal boundary curve.

\begin{figure}
  \centering
  \includegraphics[width=\textwidth]{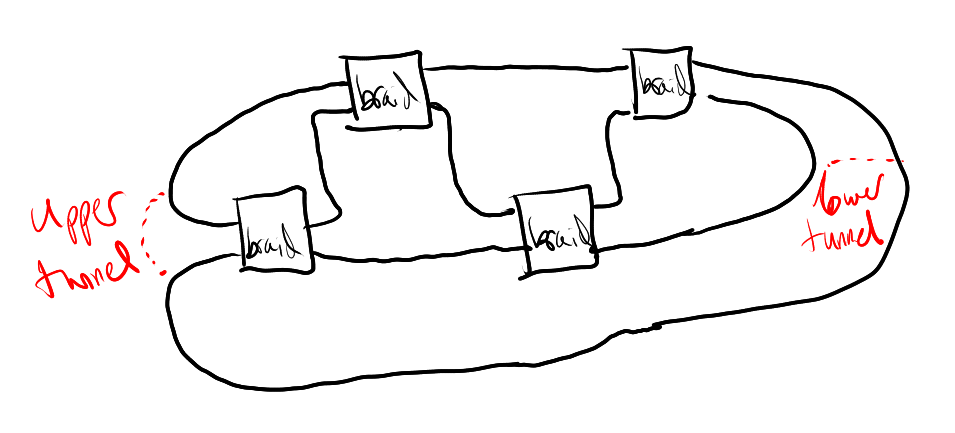}
  \caption{A 2-bridge knot with four braid regions (in black) together with upper and lower unknotting tunnels (in red).\label{fig:2bridgeunknot}}
\end{figure}
Let $ \Gamma_{\rho_0} = \langle X, Y_{\rho_0} : WX = Y_{\rho_0} W \rangle $ be a faithful representation of the 2-bridge link group of slope $ p/q $ (where $ W $ is a word in $ X $ and $ Y_{\rho_0} $ depending on the
slope; compare (\ref{eqn:overpres}) and \cite[Proposition 1]{riley72}). An \textit{unknotting tunnel} for a link $ L $ is a properly embedded arc $ t $ in the complement manifold $ M $ for $ L $
such that if $ N $ is a tubular neighbourhood of $ t $ then $ M \setminus N $ is a handlebody; a 2-bridge link has two natural unknotting tunnels, depicted in Figure~\ref{fig:2bridgeunknot}, called
the \textit{upper} and \textit{lower tunnels} \cite{sakuma98}. The element in $ \Gamma_{\rho_0} $ corresponding to a loop around the upper tunnel $ t $ is $ WXY_{\rho_0}^{-1}W^{-1} $; in $ \Gamma_{\rho_0} $,
this element is the identity. The discrete groups on the extended pleating ray which contains $ \Gamma_{\rho_0} $ correspond to the values of $ \rho \in \CC $ for which the corresponding word in $ \Gamma_\rho $ becomes
elliptic of ever-increasing order; in the limit (when the curve reaches the boundary) this word becomes parabolic and the corresponding manifold has a third deleted arc (and the Riemann surface of
the group has become a pair of disjoint 3-punctured spheres). The limit group is called a \textit{cusp group}. Pushing further along this curve, the group becomes free and the quotient surface becomes
a 4-times punctured sphere. (A more in-depth look at this process may be found in \cite[\S 2.2]{ems22}.) An alternative way of seeing that the quotient, for $ \rho $ large, is a 4-times punctured sphere is
via the \textit{Poincar\'e polyhedron theorem} (for which see \cite[\S 11.2]{ratcliffe} or \cite[IV.H]{maskit}); we explain this in a little more detail in \cite[\S 5.1]{ems22}.

\begin{figure}
  \centering
  \includegraphics[width=\textwidth]{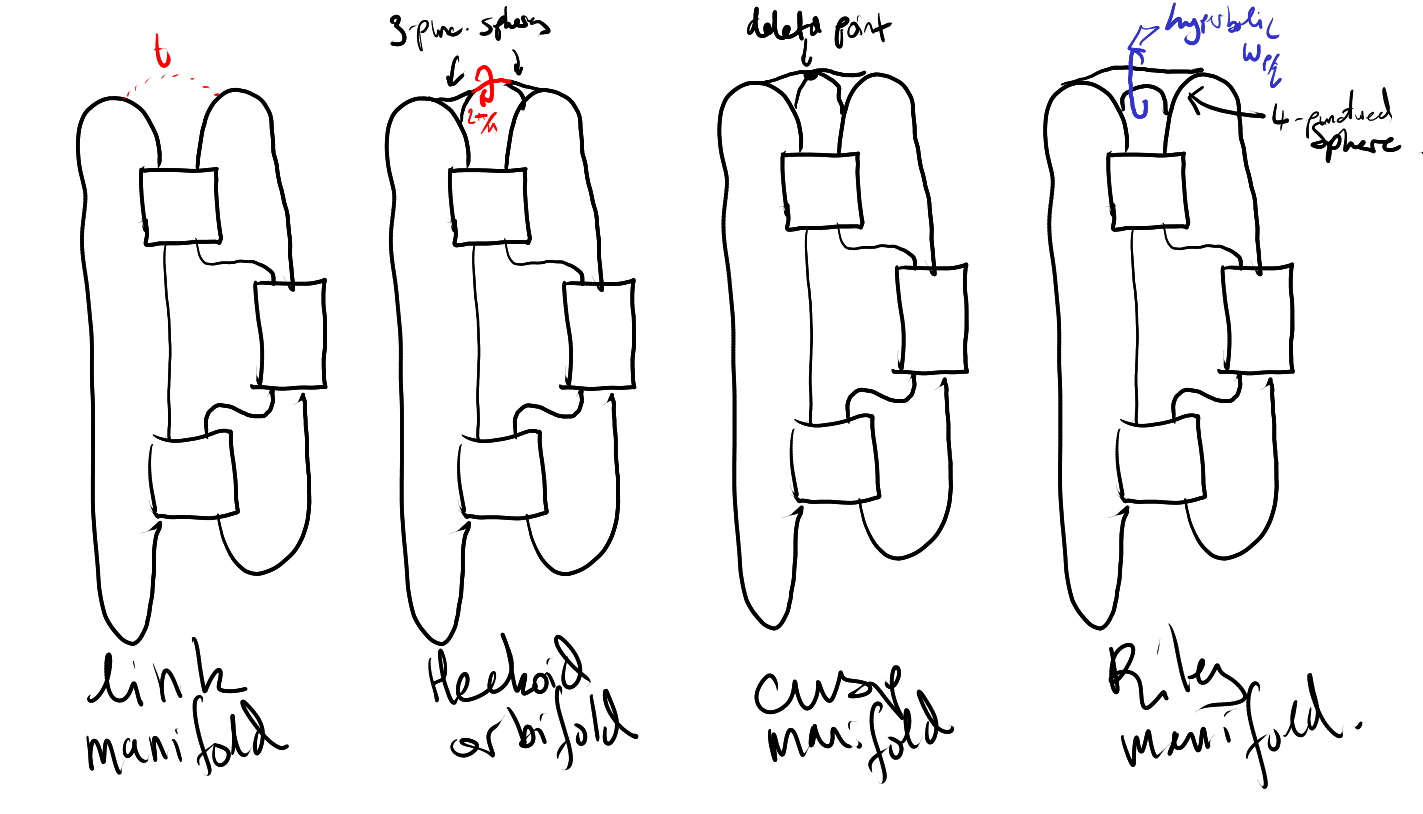}
  \caption{The four kinds of orbifolds found along an extended pleating ray.\label{fig:2bridgelimit}}
\end{figure}

The quotient manifolds which occur along the curve are sketched in Figure~\ref{fig:2bridgelimit}.
Looking back at Riley's plot, Figure~\ref{fig:riley_orig_pic}, the fractal boundary is the set of cusp groups, and the empty space in the upper right half of the plot is the moduli space of four-times punctured
spheres.

\begin{definition}\label{defn:riley_para}
  The \textit{parabolic Riley slice}, which we will often shorten to \textit{Riley slice} and denote by $ \mathcal{R} $, is the set
  of $ \rho \in \CC $ such that $ \Gamma_\rho $ is free, discrete, and the Riemann surface $ \Omega(\Gamma_\rho)/\Gamma_\rho $ is homeomorphic to a $4$-times punctured sphere.
\end{definition}
Indeed, $ \overline{\mc{R}} $ as we defined it in Section~\ref{sec:rileybounds} is the closure in $ \CC $ of this set, and $ \mc{R} $ is the interior of $ \overline{\mc{R}} $ in $ \CC $;
this is a 2010 result of Ohshika and Miyahi \cite{ohshika10} (see also Section~3.3 of \cite{ems22} for a thorough discussion of the topology of the boundary).

Finally, we remark that of course not every group $ \Gamma_\rho $ for $\rho \in \mc{R} $ is obtained by twisting out from a 2-bridge link group; the ones that can be obtained in this way
are the groups lying on the so-called \emph{rational pleating rays} of Keen and Series \cite{keen94} (see Section~\ref{sec:ks} below). The remainder of the groups in $ \mc{R} $
are obtained by pushing out from geometrically infinite groups, and can be traced back to the fundamental groups of \emph{wild} knots, c.f. \cite{kauffman04}.

\section{Quasiconformal deformation spaces}\label{sec:qc}
In this section we define formally what we mean by a moduli space of Kleinian groups; the first two sections (\ref{sec:riemann} and \ref{sec:hm}) give a fast
introduction to the general theory, and the third section (\ref{sec:rileytop}) explains what this theory means for the Riley slice.

\subsection{Riemann moduli theory}\label{sec:riemann}
We begin by quickly recapping the classical theory of Riemann moduli space; some nice introductory books on this subject include those by Farb and Margalit \cite{farb}, by
Imayoshi and Taniguchi \cite{imayoshi_taniguchi}, and by Lehto \cite{lehto}. These books all assume comfort with Riemann surfaces; an introductory text on these which has the
advantage of being cowritten by one of the founders of the moduli theory of Kleinian groups is the book by Farkas and Kra \cite{farkas_kra}. Historical accounts of
Teichm\"uller and moduli theory include \cite{acampo16,ji13,ahlfors53}; some further history may be found in \cite[pp.442--450]{segal03}.

For the Riley slice, we are only interested in finite hyperbolic Riemann surfaces (i.e. Riemann surfaces $ S $ with Euler characteristic $ \chi(S) < 0 $ and which may be compactified by filling in
finitely many punctures) and so we restrict ourselves to this case here. A discussion about what changes need to be made in order to deal
with the case of possibly non-hyperbolic Riemann surfaces with more general boundary components can be found at the start of Section~5.1.1 of \cite{imayoshi_taniguchi}.

A \textit{hyperbolic structure} on a Riemann surface $ S $ is a pair $ (X,\phi) $ where $ X $ is a surface
with a complete hyperbolic metric of finite area with totally geodesic boundary and where $ \phi : S \to X $ is a diffeomorphism. We say
that two hyperbolic structures $ (X_1,\phi_1) $, $ (X_2,\phi_2) $ are \textit{equivalent} if there exists an isometry $ \iota : X_1\to X_2 $
such that the diagram
\begin{displaymath}
\begin{tikzcd}
  &S\arrow[dl,"\phi_1"']\arrow[dr,"\phi_2"]\\
  X_1\arrow[rr,"\iota"] & & X_2
\end{tikzcd}
\end{displaymath}
commutes up to homotopy. The Teichm\"uller space of $ S $, $ \Teich(S) $, is the set of all hyperbolic structures on $ S $ modulo equivalence.
It has a natural metric space structure, but defining the metric requires quasiconformal machinery and so we put off the explanation until
Section~\ref{sec:hm}.

The \textit{mapping class group} of a Riemann surface $ S $, denoted $ \MCG(S) $ (also called the \textit{modular group} of $ S $), is the group
of orientation-preserving homeomorphisms of $ S $, moduli isotopies. (Observe that these homeomorphisms are allowed to permute the punctures. In
the general case, we define the mapping class group in such a way that punctures are permuted, but all other boundary components are fixed pointwise;
see \cite[\S 2.1]{farb} and \cite[\S 6.3.1]{imayoshi_taniguchi}.) This group has a natural action on $ \Teich(S) $: if $ f $ is a representative of
an element of the mapping class group, and $ (X,\phi) $ is a representative of an element of the Teichm\"uller space, then
\begin{displaymath}
  f \cdot (X,\phi) := (X,\phi \circ f^{-1}).
\end{displaymath}
(The use of $ f^{-1} $ rather than $ f $ is necessary for the axioms for a group action to hold.)

The mapping class group acts discontinuously as a group of isometries on $ \Teich(S) $ \cite[Theorem 6.18]{imayoshi_taniguchi}, and the quotient
$ \Teich(S)/\MCG(S) $ is called the \textit{(Riemann) moduli space} and denoted by $ \Mod(S) $. The Riemann moduli space parameterises hyperbolic
structures but now doesn't care about the position of punctures in the metric; quoting Thurston \cite[p.259]{thurstonB}:
\begin{quote}
  Informally, in Teichm\"uller space, we pay attention not just to what metric a surface is wearing, but also to how it is worn. In moduli space,
  all surfaces wearing the same metric are equivalent. The importance of the distinction will be clear to anybody who, after putting a pajama
  suit on an infant, has found one leg to be twisted.
\end{quote}

We give two examples.
\begin{example}[Tori]
  The canonical first example of a moduli space is that of genus 1 Riemann surfaces (i.e. elliptic curves). A complex structure on the torus is determined
  by the lengths of the meridian and longitude; equivalently, every complex torus is of the form $ \CC/G_\tau $ where $ \tau \in \HH^2 $ and
  \begin{displaymath}
    G_\tau := \{m + n\tau : m,n \in \ZZ\}.
  \end{displaymath}
  Further, the tori corresponding to $ \tau_1 $ and $ \tau_2 $ are biholomorphic iff $ \tau' = A\tau $ for some $ A \in \PSL(2,\ZZ) $ (acting on the upper half-plane
  as a fractional linear transformation). We therefore get the following for a complex torus $ T $:
  \begin{gather*}
    \Teich(T) = \HH^2, \quad \MCG(T) = \PSL(2,\ZZ),\\\text{and } \Mod(T) = \HH^2/\PSL(2,\ZZ) \simeq_{\text{bihol.}} S_{(2,3,\infty)}.
  \end{gather*}
  Here the space $S_{(2,3,\infty)}$ is topologically the complex plane. However the natural geometric structure it inherits as a quotient by an isometric action is hyperbolic.  $S_{(2,3,\infty)}$ is the finite area (in fact area $\frac{\pi}{3}$) Riemann surface obtained as the two sphere with two cone points of order $2$ and $3$ and a puncture.   There is a unique hyperbolic metric on this space.
  
  For further detail see \cite[\S 1.2]{imayoshi_taniguchi}, \cite[\S 10.2, \S 12.2]{farb}. As we alluded to, this is very closely related to the arithmetic study
  of elliptic curves and the theory of modular forms: see for instance \cite[C.12]{silverman}.
\end{example}

\begin{example}[Four-punctured spheres]\label{ex:fourpuncsph}
  The Teichm\"uller theory of the 4-times punctured sphere $S_{0,4}$ is intimately related to that of the once punctured torus $T_*$, since there is a natural involution (the \textit{hyperelliptic involution})
  which swaps the homotopy classes of non-boundary-parallel simple closed curves on $ T_* $ and $ S_{0,4} $ (see \cite[\S 2.2.5]{farb} and \cite[Chapter 2]{akiyoshi}).  In fact $T_*$ can be obtained much as in the previous example,
  \begin{equation}
  T_* =   \big(\CC\setminus \{G_\tau(0)\}\big)/G_\tau, \quad\quad \tau \in \HH^2  
  \end{equation}
  and if we put the complete hyperbolic metric on $\CC\setminus \{G_\tau(0)\}$,  the action of $G_\tau$ will be by isometry.  This construction gives all the complex structures on $T_*$.   In any case, it can be shown  that
  \begin{displaymath}
    \Teich(S_{0,4}) = \HH^2 \quad\text{and}\quad \MCG(S_{0,4}) = \PSL(2,\ZZ) \times (\ZZ/2\ZZ \times \ZZ/2\ZZ);
  \end{displaymath}
  one can also show that $ \Mod(S_{0,4}) $ is topologically a thrice-punctured sphere.
\end{example}

Moduli spaces of Riemann surfaces appear naturally in many interesting areas of mathematics,  through to the study of string interactions \cite[\S 22.6--22.7]{zwiebach} in physics.  We recommend the survey article by Alex Wright of Mirzakhani's work on the Moduli spaces of Riemann surfaces for an accessible update on the most  recent advances,  ideas and conjectures, \cite{AW}.

\subsection{Holomorphic motions and quasiconformal deformations}\label{sec:hm}
The analogue of the Riemann moduli space for hyperbolic 3-manifolds (moving from parameterising the conformal boundary to the interior) requires
the notion of quasiconformal mappings. The goal of this section is to give an intuitive reason for the occurence of this class of deformations
in particular, via the theory of holomorphic motions; a more detailed exposition of this route may be found in Chapter~12 of the book of Astala,
Iwaniec, and Martin \cite{astala}; it is known, following Sullivan and McMullen \cite{sullivan1,sullivan2,sullivan3,mcmullenRFC}, that the entire
theory is basically analogous to the theory of complex dynamics and so the reader should try to focus on the deformations occuring in $ \hat\CC $
(in this case, deformations of limit sets) as well as deformations of the hyperbolic 3-manifolds: we find it much easier anyway to visualise points moving around in $ \hat\CC $ compared with
trying to think about deforming large hyperbolic manifolds. The idea is that, although the Riemann surface cannot detect all of the Kleinian structure (the Riemann moduli space is a quotient
of the Kleinian moduli space which we will define), the limit set \emph{can} (as it should, since it is a global reflection of the group orbits).

\begin{figure}
  \centering
  \includegraphics[width=\textwidth]{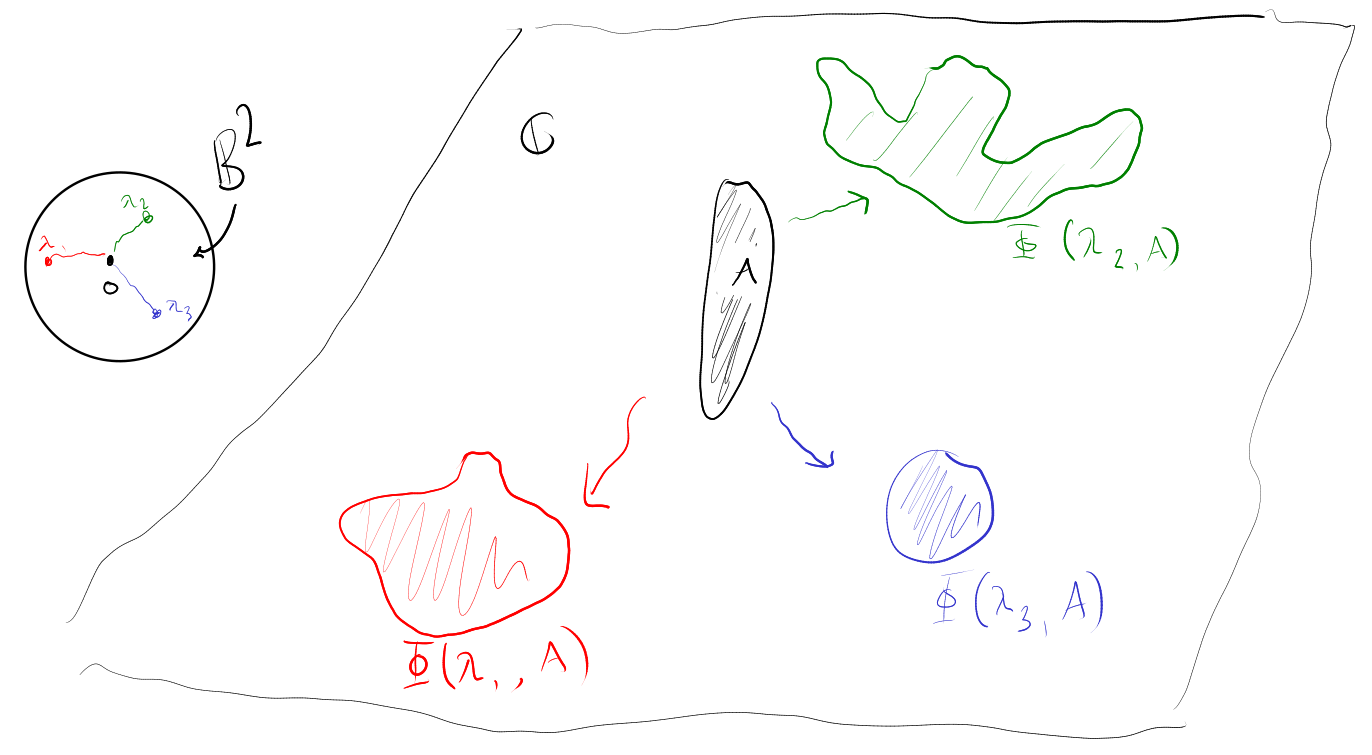}
  \caption{Images of $ A \subset \CC $ under a holomorphic motion $ \Phi $.\label{fig:holomotion}}
\end{figure}

The basic definition is due to Mañé, Sad, and Sullivan \cite{mane83}; it encodes what we mean when we say that a parameter `varies holomorphically' in a set.
\begin{definition}[Holomorphic motions]
  Let $ A \subseteq \hat\CC $. A \textbf{holomorphic motion} of $ A $ is a map $ \Phi : \BB^2 \times A \to \hat\CC $ (where $ \BB^2 $ is the unit disc in $ \CC $) such that
  \begin{enumerate}
    \item For each $ a \in A $, the map $ \BB^2 \ni \lambda \mapsto \Phi(\lambda,a) \in \hat\CC $ is holomorphic;
    \item For each $ \lambda \in \BB^2 $, the map $ A \ni a \mapsto \Phi(\lambda,a) \in \hat\CC $ is injective;
    \item The mapping $ A \ni a \mapsto \Phi(0,a) \in \hat\CC $ is the identity on $ A $.
  \end{enumerate}
  See the schematic in Fig.~\ref{fig:holomotion}.
\end{definition}

Mañé, Sad, and Sullivan also proved the so-called `$\lambda$-lemma'; we need an extended version (Theorem~\ref{thm:extended_lambda_lemma}), due
to S\l{}odkowski. The result shows that holomorphic motions are `rigid': they are determined everywhere on $ \hat\CC $ even when only known
on `small' sets, in such a way that the extension is (remarkably) almost conformal in $ a $. More precisely, its action on the tangent bundle of $ \CC $ maps circles
to ellipses of bounded eccentricity---this is the germ of the notion of quasiconformality, which we pause to define following Section~1.4 and Chapter~4
of \cite{imayoshi_taniguchi}. An extended treatment can be found in \cite{astala}, but this is perhaps too technical for bedtime reading.

\begin{definition}[Quasiconformality]
  Let $ D, D' $ be domains in $ \CC $, and let $ f : D \to D' $ be an orientation-preserving diffeomorphism. Define the two differential operators
  \begin{displaymath}
    \partial_z := (1/2)(\partial_x - i\partial_y) \quad\text{and}\quad \partial_{\overline{z}} := (1/2)(\partial_x + i\partial_y)
  \end{displaymath}
  where $ x = \Re z$ and $ y = \Im z $. The function
  \begin{displaymath}
    \mu = \frac{\partial_{\overline{z}} f}{\partial_z f}
  \end{displaymath}
  is a smooth complex-valued function on $ D $, and $ \|\mu\|_\infty < 1 $. The function $ \mu $ is called the \textit{Beltrami coefficient} of $ f $, and
  the function $ f $ is called a \textit{quasiconformal mapping} if
  \begin{displaymath}
    K(f) := \sup_{z\in D} \frac{1+|\mu(z)|}{1-|\mu(z)|} < \infty
  \end{displaymath}
  in which case $ K(f) $ is called the \textit{maximal dilatation} of $ f $.
\end{definition}

\begin{figure}
  \centering
  \includegraphics[width=0.6\textwidth]{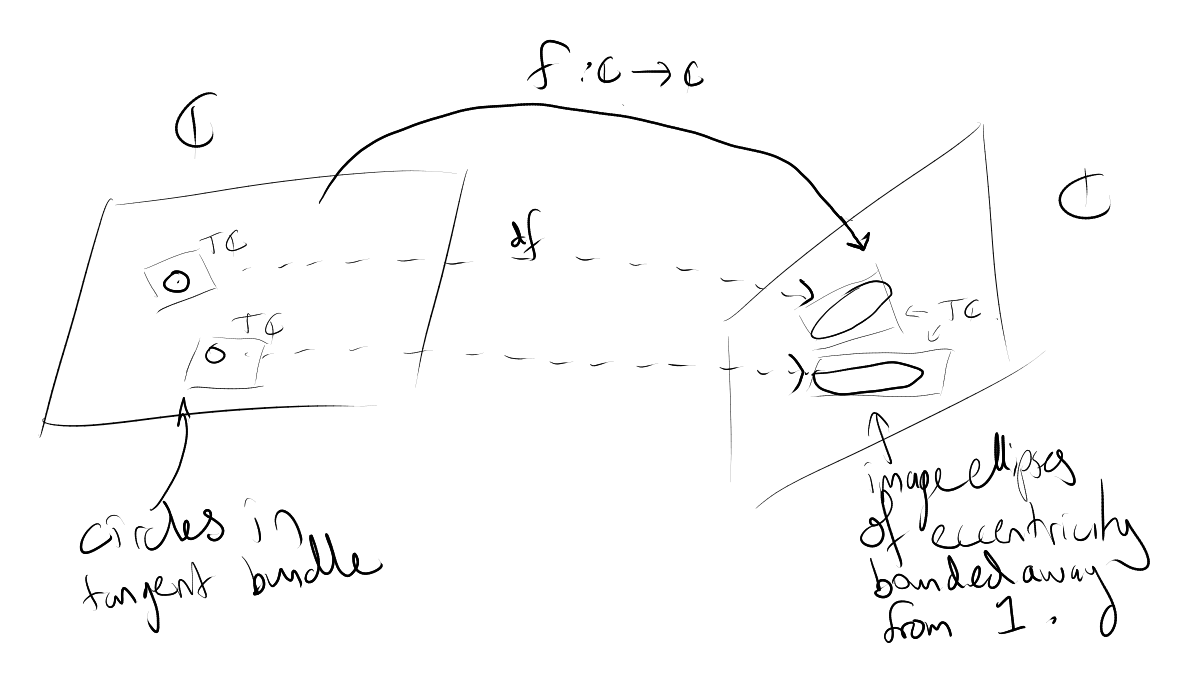}
  \caption{The image of circles in the tangent bundle to $ \CC $ under the differential of a quasiconformal map $ f : \CC \to \CC $.\label{fig:quasiconformality}.}
\end{figure}

This definition relates to the action on tangent bundle circles in the sense that $ \mu(z) $ measures the eccentricity of ellipses which are the images
of circles in $ T_z \CC $ under the differential $ df $ and $ K(f) $ measures the maximal ratio of one axis of these ellipses to the other,
see \cite[p.18]{imayoshi_taniguchi} and Fig.~\ref{fig:quasiconformality}.

It turns out that this notion is the `correct' weakening of conformality to enable the extension of holomorphic motions.

\begin{theorem}[{S\l{}odkowski's extended $ \lambda$-lemma \cite{Slod1,Slod2}}]\label{thm:extended_lambda_lemma}
  If $ \Phi : \BB^2 \times A \to \hat\CC $ is a holomorphic motion of $ A \subseteq \CC $, then $ \Phi $ has an extension to $ \tilde \Phi : \BB^2 \times \hat\CC \to \hat\CC $
  such that
  \begin{enumerate}
    \item $ \tilde\Phi $ is a holomorphic motion of $ \hat\CC $;
    \item For each $ \lambda \in \BB^2 $, the map $ \tilde\Phi_\lambda $ defined by $ \hat\CC \ni a \mapsto \tilde\Phi(\lambda,a) \in \hat\CC $ is a quasiconformal homeomorphism with
          maximal dilatation satisfying
          \begin{displaymath}
            K(\tilde\Phi_\lambda) \leq \frac{1+|\lambda|}{1-|\lambda|};
          \end{displaymath}
    \item $ \tilde\Phi $ is jointly continuous in $ \BB^2 \times \hat\CC $; and
    \item For all $ \lambda_1, \lambda_2 \in \BB^2 $, $ \tilde\Phi_{\lambda_1} \tilde\Phi_{\lambda_2}^{-1} $ is quasiconformal with
          \begin{displaymath}
            \log K\!\!\left(\tilde\Phi_{\lambda_1} \tilde\Phi_{\lambda_2}^{-1}\right) \leq \rho(\lambda_1,\lambda_2)
          \end{displaymath}
          (where $ \rho $ is the hyperbolic metric on $ \BB^2 $). \qed
  \end{enumerate}
\end{theorem}

There is also an equivarant version due to Earle, Kra, and Krushkal', appearing as Theorem~1 of their paper \cite{EKK}:
\begin{theorem}[Equivariant $ \lambda$-lemma]
  Let $ A \subseteq \hat\CC $ have at least three points, and let $ \Gamma $ be a group of conformal motions preserving $ A $.
  Let $ \Phi : \BB^2 \times A \to \hat\CC $ be a holomorphic motion on $ A $, and suppose that for each $ \gamma \in \Gamma $ and each $ \lambda \in \BB^2 $
  there is a conformal map $ \theta_\lambda(\gamma) $ such that
  \begin{equation}\label{eqn:equiv_hm}
    \Phi(\lambda,\gamma(z)) = \theta_\lambda(\gamma) (\Phi(\lambda,z))
  \end{equation}
  for all $ z \in A $. Then $ \Phi $ can be extended to a holomorphic motion on $ \hat\CC $ such that (\ref{eqn:equiv_hm}) holds for all $ z \in \hat\CC $. \qed
\end{theorem}

Observe that, by the $\lambda$-lemmas, a holomorphic motion of the limit set of a Kleinian group extends quasiconformally to the entire Riemann sphere; in particular,
the ordinary set. This motivates:
\begin{definition}
  The \textit{quasiconformal deformation space} of a Kleinian group $ \Gamma $, denoted by $ \QH(\Gamma) $,
  is the space of representations $ \theta : \Gamma \to \PSL(2,\CC) $ (up to conjugacy) such that
  \begin{enumerate}
    \item $ \theta $ is faithful and $ \theta\Gamma $ is discrete;
    \item $ \theta $ is type-preserving, that is if $ \gamma \in \Gamma $ is parabolic (resp. elliptic of order $ n $) then $ \theta\gamma $ is parabolic (resp. elliptic of order $ n $); and
    \item the groups $ \theta\Gamma $ are all quasiconformally conjugate (i.e. there exists some quasiconformal $ f : \hat\CC \to \hat\CC $ depending on $ \theta $ such that, as functions, $ \theta\Gamma = f\Gamma f^{-1} $).
  \end{enumerate}
  The space is equipped with a natural metric defined in roughly the same way as the classical Teichm\"uller metric (the distance between two deformations is defined to be the log of the
  maximal dilatation of the composition of the two quasiconformal homeomorphisms).
\end{definition}
This definition was studied first by Bers \cite{bers70, bers72}, Kra \cite{kra72,kra74} and Maskit \cite{maskit70,maskit71}; modern textbooks and monographs which discuss this
theory include \cite{matsuzakitaniguchi} (essentially the whole book), \cite[Chapter 5]{marden}, and \cite[Chapter 8]{kapovich}

The exact relationship between $ \QH(\Gamma) $ and $ \Teich(\Omega(\Gamma)/\Gamma) $ is subtle to define (see the references following the theorem), so we give only a rough version
here. This result is attibuted by Bers \cite[\S2.4]{bers72} to Bers and Greenberg \cite{bers71} and Marden \cite{marden71}.
\begin{theorem}\label{thm:qcds}
  Let $ \Gamma = \langle \gamma_1,...,\gamma_n \rangle $ be a finitely generated non-elementary\footnote{An `elementary' group is a Kleinian group $G$ which is particularly simple, in the sense
  that $ |\Lambda(G)| < \infty $; basically this happens when most of the fixed points of group elements collide. We do not need to worry so much since it turns out that $ \Gamma_\rho $
  is elementary iff $ \rho = 0 $, and this is excluded from the Riley slice by Lemma~\ref{lem:crudebound}. For more on elementary groups see \cite[Chapter V]{maskit} or \cite[\S 5.5]{ratcliffe}.}
  Kleinian group with $ \Omega(\Gamma) \neq \emptyset $; let $ S = \Omega(\Gamma)/\Gamma $. Then there is a well-defined holomorphic surjection $ p : \Teich(S) \to \QH(\Gamma) $. Further,
  there is a discrete subgroup $ \widehat{\MCG}(S) \leq \MCG(S) $ and a natural bijection $ \QH(\Gamma) \approx \Teich(S)/\widehat{\MCG}(S) $
  such that the two projection maps agree. \qed
\end{theorem}
In the case that $ \Gamma $ is \textit{geometrically finite} (i.e. its quotient manifold
arises as the projection of a tiling of $ \HH^3 $ by finite-sided hyperbolic polyhedra \cite[VI]{maskit}) , the group $ \widehat{\MCG}(S) $ is the subgroup generated by Dehn twists along simple closed curves which bound
compression discs.   

\medskip

The torsion-free version of this theorem appears in \cite{matsuzakitaniguchi} as Theorem 5.27 and the following discussion. The proof given there actually works for the torsion case as well, as long as we are careful to replace all manifold theorems and definitions with the corresponding `orbi-theorems' and `orbi-definitions'. Such a version
of this theorem is found as Theorem~5.1.3 of \cite{marden} (the proof is sketched as Exercise~5-35, p.367ff).

\medskip

In the case that a Kleinian group $ G $ is \textit{Fuchsian} (that is it leaves a round disc $\Delta$ in $ \Omega(G) $ invariant---and hence acts as a group of hyperbolic isometries on $ \Delta $,
with $ \Lambda(G) \subseteq \partial \Delta $, see the latter chapters of \cite{beardon} and \cite[V.G, IX, and X]{maskit}) then the quasiconformal deformation space $ \QH(G) $ is given
the special name of \textit{quasi-Fuchsian deformation space}; a Kleinian group $ \Gamma $ is called \textit{quasi-Fuchsian} if it is a quasiconformal conjugate of some Fuchsian group.
The limit sets of quasi-Fuchsian groups are contained in quasiconformal images of circles (called \textit{quasicircles}) and the general theory of these deformation spaces is fairly
well-understood: for instance, see \cite[Chapter 4]{marden} and the (very accessible) \cite{indras_pearls}. The limit set of Fig.~\ref{fig:indra} is the limit set of a quasi-Fuchsian group,
and shows how intricate these quasicircles can become; controling their large-scale geometry is one of the themes of our preprint \cite{ems21}, in which we study how small quasi-deformations
of limit sets modify the quotient surface (see also Section~\ref{sec:nbds} below).

\subsection{The Riley slice as a quasiconformal deformation space}\label{sec:rileytop}
In Section~\ref{sec:freegps} and Section~\ref{sec:riley}, we introduced the Riley slice $ \mc{R} $ as a subset of the complex plane $ \CC $ with no additional structure. However, this
set was defined in a way that seems at first glance to be compatible with holomorphic motions: if $ \Gamma_\rho $ is free and discrete with quotient surface a 4-times punctured sphere, then
surely so is $ \Gamma_{\tilde\rho} $ for a small holomorphic motion $ \tilde\rho $ of $ \rho $? That this is the case is actually surprisingly nontrivial to prove; we give full details
in our preprint \cite{ems22}, and in this section we give some of the history of the study of the Riley slice from a quasi-point of view.

The first mention that $ \mc{R} $ is a quasiconformal deformation space which we can find in the literature is in the proof of Theorem~1 of Lyubich and Suvorov's 1988 paper \cite{lyubich88}, though the proof is only
sketched; a more detailed proof along the same lines can be found in Section~3.3 of our preprint \cite{ems22}. Using this equivalence, they prove that:
\begin{itemize}
  \item $ \mc{R} $ is connected and has connected complement---note that their $ \Gamma $ is our $ \CC\setminus\mc{R} $ \cite[Theorem~1]{lyubich88};
  \item There are only two maximal parabolic conjugacy classes in $ \Gamma_\rho $ for $ \rho \in \mc{R} $ \cite[Lemma~1]{lyubich88};
  \item The Riley slice is invariant under a particular dynamical system \cite[Lemma~2]{lyubich88}, c.f. \cite{martin19}, \cite{martin20}, \cite{martin21}, \cite[Section~3]{ems21}, and our upcoming preprint \cite{ems22b};
  \item The set of $ \rho $ such that $ \Gamma_\rho $ is non-discrete is dense in $ \CC\setminus\mc{R} $ \cite[Theorem~2]{lyubich88}.
\end{itemize}

The same proof (again, only sketched) was given by Maskit and Swarup in 1989 \cite[\S 1]{maskit89}; in this paper, they show that if $ \rho \in \overline{\mc{R}} $
and $ \Omega(\Gamma_\rho) \neq \emptyset $ then the Riemann surface $ \Omega(\Gamma_\rho)/\Gamma_\rho $ is either a 4-times punctured sphere, or a disjoint union of 3-punctured spheres.
Thus, since quasiconformally conjugate groups have homeomorphic quotient surfaces, the boundary $ \partial \mc{R} $ contains either groups with empty ordinary set, or groups with
surface a pair of 3-punctured spheres. The groups in the latter class are called the \textit{cusp groups}; they are exactly the \textit{maximal cusps} of the Bers deformation
theory \cite{bers70,maskit70,maskit83,ohshika98}. The 1989 paper also shows that every group in $ \mc{R} $ is geometrically finite; thus we can apply the theorem
relating the quasiconformal deformation space with the Teichmüller space (Theorem~\ref{thm:qcds}), along with Example~\ref{ex:fourpuncsph}, to see that $\mc{R}$ is homeomorphic
to an (open) annulus---this follows from the observation that $ \widehat{\MCG}(\Gamma_\rho) $ is generated by a single parabolic element \cite[Corollary~3.9]{ems22}.

Finally we note that the Riley slice is not just a quasiconformal deformation space, it is a \emph{quasi-Fuchsian} deformation space: when $ \rho \in \RR\cap\mc{R} $, the group $ \Gamma_\rho $
is Fuchsian. (A detailed discussion of this point may be found in Section~6.3.1 of \cite{ems22}.)

\section{A geometric coordinate system}\label{sec:ks}
In this section we motivate and describe the Bers-Thurston deformation theory of quasi-Fuchsian groups in the geometrically finite case (\ref{sec:laminations}),
and then describe how the Riley slice is a concrete example of this theory via the results of Keen and Series (\ref{sec:coordinatesystem}).

\subsection{Measured laminations and cusps}\label{sec:laminations}
Given a quasi-Fuchsian space $ \mc{F} $, there are two natural definitions for the boundary $\partial\mc{F} $.
\begin{enumerate}
  \item The \textit{Bers boundary}, which is defined analytically as the boundary of a particular embedding of $ \mc{F} $ into a particular Banach space of $\Gamma$-automorphic
        forms (originally defined by Bers \cite{bers70}, see also \cite[Section~4.2]{matsuzakitaniguchi} and \cite[\S 5.11.1]{marden}); and
  \item the \textit{Thurston boundary}, which is defined geometrically in terms of measured laminations; this appears in \cite{thurston98}, see also \cite[Chapter~6]{matsuzakitaniguchi}
        and \cite[\S 5.11.2]{marden}.
\end{enumerate}

In this paper, we are interested in the Thurston boundary; and in this section we give a vague intuitive description which will be sufficient for motivating the Keen--Series theory in the
next subsection. For motivation of our description of this boundary, we will use one of the major conjectures of hyperbolic 3-manifold theory posed by Thurston \cite[Problem 11 in \S 6]{thurston82} as a very strong generalisation of Mostow rigidity;
namely, the \textit{ending lamination conjecture/theorem}. This conjecture was ultimately proved by Brock, Canary and Minsky \cite{elt2}, and see earlier work of Minsky \cite{elt1}, in the tame case; that this is sufficient follows from the
resolution of Marden's tameness conjecture (posed in \cite{marden74}, proved independently by Agol \cite{agol04} and Calegari and Gabai \cite{calegari04}).

Our exposition follows \cite{marden} (though see also \cite[Chapter 14]{kapovich}, \cite[Section 4]{series05}, and \cite{thurstonN}); we will explain all of the words used in the theorem statement after giving it.
\begin{theorem}[Ending lamination theorem]\label{thm:elt}
  A hyperbolic 3-manifold with finitely generated fundamental group is uniquely determined by its topological type and its end invariants.
\end{theorem}
Here, `topological type' just means homeomorphism class. Thus the only thing we need to explain is what we mean by `end invariant'. This consists
of two things: an \textit{end} of the manifold, and a \textit{measured lamination} on that end. An end of a 3-manifold, roughly speaking, is an ideal boundary
component of the manifold; in the Riley slice interior (the main case which we are interested in), all groups are geometrically finite and in this case all ideal ends actually `exist' in the sense
that they are just boundary components in the usual sense (the components of the quotient Riemann surface).

\begin{figure}
  \centering
  \includegraphics[width=0.5\textwidth]{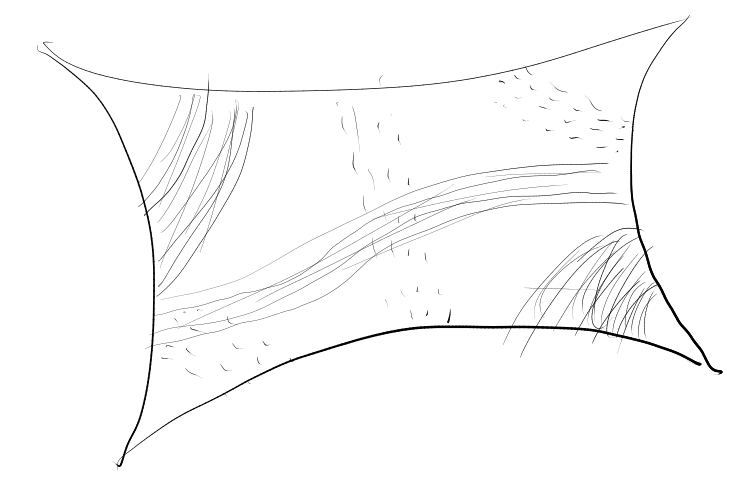}
  \caption{A lamination of infinite length on the 4-times punctured sphere.\label{fig:inflam}}
\end{figure}

The definition of a measured lamination requires more work. A \textit{geodesic lamination}---or simply \textit{lamination}---on a surface $ S $ is a compact subset
of $ \interior S $ which is a union of complete disjoint geodesics (each geodesic is called a \textit{leaf} of the lamination); an example of a lamination with infinitely many windings around the surface
is given in Fig.~\ref{fig:inflam}, and many other pictures are found in \cite[Chapter 6]{matsuzakitaniguchi}. A \textit{measured lamination} is a pair $ (\Lambda,\mu) $ where $ \Lambda $ is a lamination,
and $ \mu $ is a \textit{transverse measure} supported on $ \Lambda $. The reader is not expected to be familiar with the ins-and-outs of measure theory; we only need
the general idea that a measure is `a density which accumulates on curves when integrating'. (A good reference, though, is Rudin's second textbook \cite{rudinRCA}.) In
this case we wish to assign to every curve $ \gamma $ which is not tangent to any leaf of $ \Lambda $ (`transverse') a number $ \int_\gamma d\mu $ in a consistent way (`measure'). The set
of measured laminations on $S$ admits a natural topology (the weak convergence topology) and the resulting topological space is denoted by $ \ML(S) $.

\begin{figure}
  \centering
  \includegraphics[width=\textwidth]{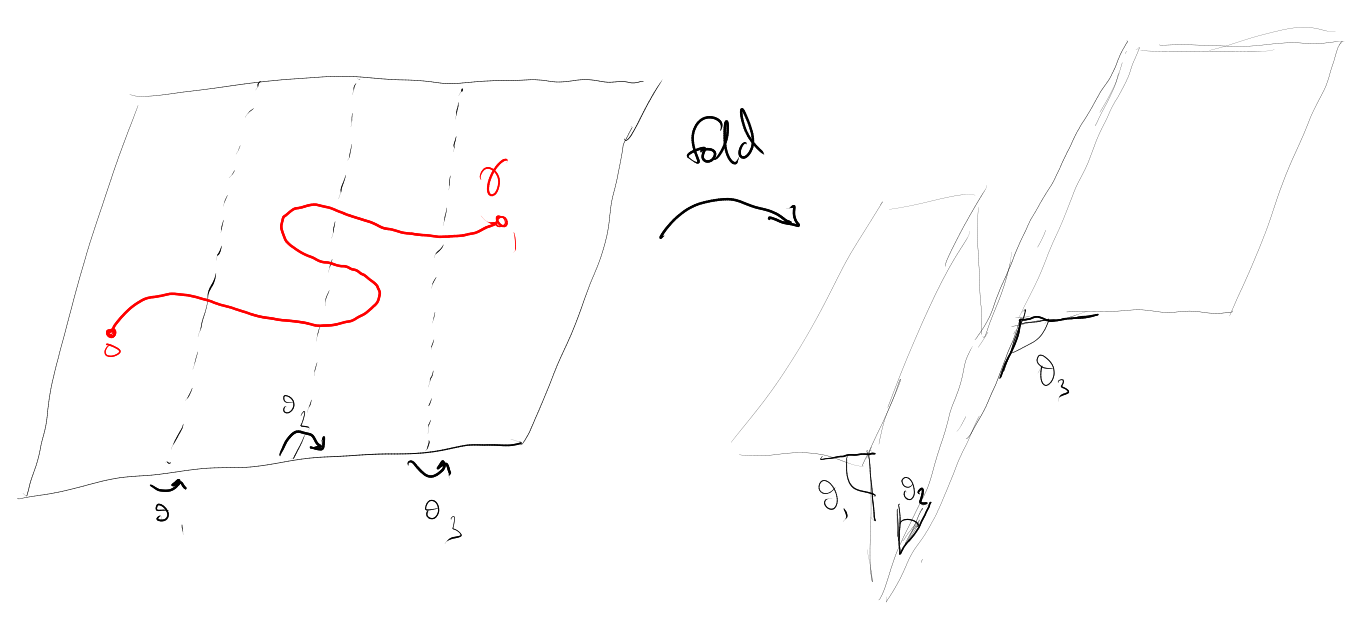}
  \caption{A (Euclidean) pleated surface, with pleating locus (set of folds) consisting of three geodesics. The red curve $ \gamma $ (oriented from 0 to 1) has
  measure $ \theta_1+\theta_2-\theta_2+\theta_2+\theta_3 = \theta_1+\theta_2+\theta_3 $ with respect to the induced measured lamination.\label{fig:pleatedsurface}}
\end{figure}

\begin{example}[Pleated surfaces]
  Take a surface $ S $ and fold it along some geodesics so that along each geodesic the folding angle is constant. If the union of the geodesics forms a lamination $ \Lambda $,
  then the \textit{pleating measure} on $ \Lambda $ is the measure which assigns to each curve on $ S $ (...which is not tangent to any of the folds...) the sum of the folding angles
  of the folds which it crosses.\footnote{That all this is well-defined is really non-trivial, see Thurston's description \cite[\S 8.6]{thurstonN} and the paper by Canary, Epstein, and Green \cite{canary}.}
  We draw a Euclidean example with three folds in Fig.~\ref{fig:pleatedsurface}. A very common way that one runs into pleated surface is by considering the hyperbolic convex hull of the limit set of a Kleinian group; two very nice visualisations of this have been produced by Brock and Dumas \cite{jewel,bug}.
\end{example}

An ending invariant of a Kleinian group $ G $ in the sense of Theorem~\ref{thm:elt} (in the geometrically finite case...) is just a connected component $ \Omega^* $ of $ \Omega(G)/G $ together
with a measured lamination on $ \Omega^* $.

We now proceed to describe the \textit{Thurston boundary} of quasi-Fuchsian space. A measured lamination $ (\Lambda,\mu) $ on a hyperbolic surface $ S $ has a natural \textit{lamination length}, defined to be
the measure on $ S $ which is locally the product of $ \mu $ (in the direction locally transverse to leaves of $ \Lambda $) and the usual hyperbolic length measure (in the direction
locally parallel to $ \Lambda $). If we change the complex structure on $ S $ by shrinking the length of $ (\Lambda,\mu) $ then (intuitively) the topological type of the surface does not
change until the length is pinched down to 0; neighbourhoods of the leaves of $ \Lambda $ will turn into neighbourhoods of deleted points on the surface (if the leaves are discrete),
and we have hit the boundary of the Teichm\"uller space of the surface. It turns out that this process gives \emph{quasi-conformal} deformations of the underlying group, and (if everything
is well-behaved) the union of all the surfaces obtained in this way forms a boundary of the deformation space (in the sense that these boundary surfaces, together with the deformation
space itself, form a compact set).

Now observe that all of the surfaces which we have produced on the boundary are essentially different (they have different complex structures); thus, we have a natural bijection between
this boundary and the set of measured-laminations-up-to-length: this latter set is usually called the space of \textit{projective classes of measured laminations}, $ \PML(S) $.

All of this (the geometrically finite ending lamination theorem and the Thurston compactification theory) gives the following `picture' of the quasiconformal deformation space $ \mc{F} $:
it consists of a set of groups representing a family of surfaces, where the surfaces are all of the same topological type and where deforming the group slightly either changes the length
of a distinguished set of geodesics on the surface (this is a `radial' deformation, towards the boundary) or changes from one distinguished set of geodesics to another which is `close'
in the sense both of the closeness of the physical set of geodesics and of the measure on the geodesics. The boundary is the set of groups which parameterise those surfaces in which
these geodesics have been pinched to deleted points. There is a subtlety here, in the sense that this procedure \emph{a priori} will only produce the geometrically finite groups on the boundary since
a measured lamination might correspond to a geodesic of infinite length, and so some normalisation needs to take place to make sense of `decreasing the length'; in the case of the
deformation space of a once-punctured torus, this is explained in \cite[\S 6]{keen93} and in the case of the Riley slice we explain the process in \cite[\S 6.3 -- 6.4]{ems22}. The
groups obtained by pinching down geodesics of finite length are called the \textit{cusp groups}, and if their corresponding surfaces are unions of 3-punctured spheres then they are
called \textit{maximal cusps}; by a result of McMullen \cite{mcmullen91}, these are dense in the boundary. In terms of the group theory, the geodesics are initially represented by
loxodromic elements (the elements arising in the group as so-called \textit{Farey words}), the pinching corresponds to deforming these elements
so their trace approaches $ -2 $, and at the cusps the elements become parabolic. This process is studied in general in \cite{bers70,maskit70,maskit83,ohshika98}. The groups produced by
pinching geodesics of \emph{infinite} length have geometrically infinite ends, and correspond to the non-cusp points on the boundary; a picture of this case is found as Fig.~6.5
of \cite{matsuzakitaniguchi} (in the case of the 4-times punctured sphere, there are no finite ends after pinching: this can be seen by looking at a generic infinite lamination, like Fig.~\ref{fig:inflam}).

\subsection{The Keen--Series coordinate system}\label{sec:coordinatesystem}
We now discuss the measured lamination theory in the specific case of the Riley slice; this was developed by Keen and Series \cite{keen94} with minor corrections
by Komori and Series \cite{komoriseries98}; some of the theory was not developed in detail in those papers as it is similar to the case of the once-punctured torus
moduli space (the \textit{Maskit slice}) which Keen and Series had studied earlier \cite{keen93,series99}. A full account of the theory may be found in our preprint \cite{ems22},
though we develop the elliptic case at the same time which adds some difficulty for the reader.

In any case, the first main theorem is the following.
\begin{theorem}[{\cite[Theorem~5.4]{keen94}}]\label{thm:ksmain}
  There is a natural homeomorphism
  \begin{displaymath}
    \Pi : \mc{R} \to \RR/2\ZZ \times \RR_{>0}.
  \end{displaymath}
\end{theorem}

To prove this theorem, one first shows that $ \PML(S_{0,4}) $ is naturally homeomorphic to $ \RR/2\ZZ $ (see e.g. Section~2 of \cite{keen94} and Section~5 of \cite{ems22}---each measured lamination on
the surface can be assigned a `sloping angle' $ \pi \theta $ and the homeomorphism sends $ \Lambda \mapsto \theta(\Lambda) $) and then defines $ \Pi $ to be the map
which sends $ \rho \in \mc{R} $ to the pair $ (\theta(\Lambda), l) $ where $ \Lambda $ is the lamination on $ \RR $ coming from the ending lamination theorem (which can be found very concretely in this
case: if $ H $ is the hyperbolic convex hull of the limit set $ \Lambda(\Gamma_\rho) $ in $ \HH^3 $, then $ H/\Gamma_\rho $ is a deformation retract of $ \HH^3/\Gamma_\rho $ whose boundary is a pleated
surface; the lamination $ \Lambda $ is the corresponding measured lamination), and $ l $ is its lamination length. The continuity then follows from various known results about continuity of functionals
associated to measured laminations, discussed for example in \cite{keenseriesCCHB}.

In the process of showing Theorem~\ref{thm:ksmain}, Keen and Series show that the fibres of $ \Pi $ in $ \mc{R} $ corresponding to fixed angles in $ \RR/2\ZZ $ are analytic curves
in $ \CC $ (the \textit{pleating rays}, corresponding to pinching a single measured lamination) and in the case that it is a lamination of finite hyperbolic length that is being pinched
(equivalently, deforming along a fibre above $ \QQ/2\ZZ $) these curves are branches of the inverse image of $ (-\infty,-2) $ under a certain polynomial in one complex variable.
We call these polynomials the \textit{Farey polynomials} since they are very closely related to the sequences in number theory named after Farey; see the paper \cite{akiyoshi2020classification},
the monograph \cite{akiyoshi}, and our preprint \cite{ems22b}.

\begin{figure}
  \centering
  \includegraphics[width=\textwidth]{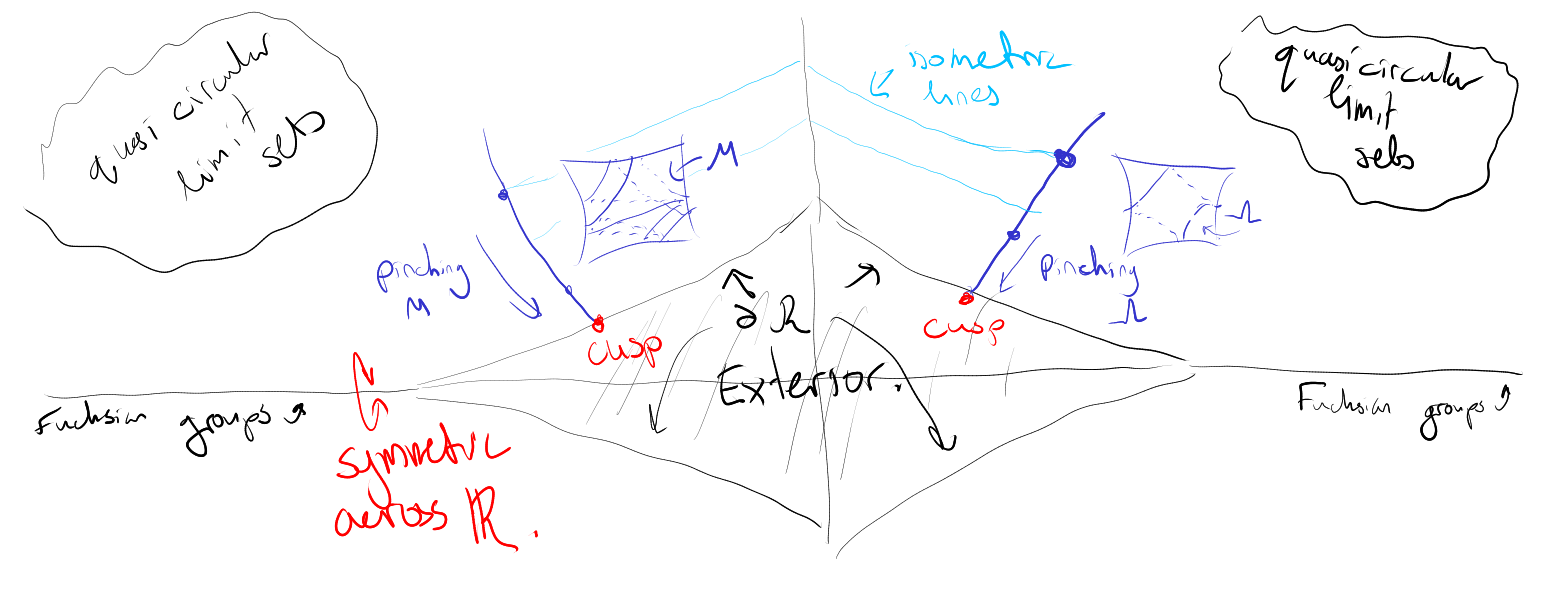}
  \caption{A map of the Keen--Series coordinate system of the Riley slice.\label{fig:rileyslicemap}}
\end{figure}

A rough map of the Riley slice may be found in Fig.~\ref{fig:rileyslicemap}. Various people have produced computer-generated images of the coordinate system
and the rational pleating rays; for example Wright \cite[Fig.~1]{gilman08} and Yamashita \cite[Fig.~0.2b]{akiyoshi}, \cite[Fig.~1]{ems21}.

We believe that the following problem is open.
\begin{problem}
  Prove or disprove that $ \Pi $ is conformal.
\end{problem}
If it is true that $ \Pi $ is conformal, then the coordinate lines of the Keen--Series coordinate system are Teichm\"uller geodesics (see \cite[\S 10.2]{thesis} for an argument).  We think the answer is probably no,  however we make a specific closely related conjecture.  Precise definitions and various facts concerning quasigeodesics can be found in  \cite{bridson_haefliger}.  Roughly $\alpha$ is a quasigeodesic if it is the image of a geodesic $\beta$ under quasiisometry. $f$, maps which admits constants $\lambda$ and $\epsilon$ so that 
\[ \frac{1}{\lambda} d(z,w)-\epsilon \leq d(f(z),f(w)) \leq {\lambda} d(z,w)+\epsilon, \quad z,w\in \beta \]
Notice such maps need not even be continuous, but they are injective and bilipschitz ``at large scales''.
 \begin{conjecture}
The rational pleating rays are quasigeodesics in the Teichm\"uller metric.
\end{conjecture}

\begin{remark}
We note that similar structures and coordinates have been worked out in detail for many other moduli spaces beyond the Maskit and Riley slices we have mentioned above; for instance, the
  so-called \textit{diagonal slice of Schottky space} \cite{series17}, the \textit{Earle embedding} \cite{komori01}, and various other
  linear slices through complex Fenchel-Nielsen coordinates for the deformation space of quasi-Fuchsian once-punctured torus groups \cite{komori12}.
\end{remark}

\section{Our recent work}\label{sec:ourwork}
In this final section we announce some of our recent work on the Riley slice; this section is primarily written for people familiar with the Riley slice theory and we rely on the previous sections
to provide historical and mathematical context.

\subsection{Extension of the Keen--Series theory to the elliptic case}\label{sec:elliptic}
There are two natural ways to motivate the study of groups generated by two elliptic elements rather than two parabolic elements.
\begin{enumerate}
  \item A group $ \Gamma_\rho $ for $ \rho \in \mc{R} $ is a \textit{Schottky type group}: there are two disjoint tangent pairs of circles in $ \hat\CC $, say $ C_1,C_1' $
        and $ C_2,C_2' $ (where $ C_i $ and $ C_i' $ are tangent for each $ i $) such that $ X(C_1) = C_1' $, $ Y_\rho(C_2) = C_2' $, and $ X $ and $ Y_\rho $ send the inside
        of the `domain circle' to the outside of the `target circle' (see e.g. \cite[Chapter 4]{indras_pearls}, \cite[Exercises IV.J.21 and VII.F.9]{maskit}, and \cite[\S 3.5]{thesis}
        for some properties of Schottky groups). These groups
        arise naturally from the study of (classical) Schottky groups, where the paired circles are assumed disjoint and the quotient is a handlebody, by pushing the circles together and hence deforming the
        corresponding groups (in fact, by some work of Maskit \cite{maskit81}, the Riley slice is a slice through the boundary of the space of Schottky groups---hence the name).
        If this deformation is continued and the circles are pushed through each other to form two pairs of transversely intersecting circles, then the pairing elements become
        elliptic.
  \item If one prefers to motivate the Riley slice via knot theory, it is natural to consider knot orbifolds rather than knot complement manifolds; more precisely,
        given a 2-bridge knot $ k $ we may produce a manifold $ M $ which is a branched cover of the knot complement where the knot outlines `portals' through which
        one may step between the different covering sheets; since $ \pi_1(S^3\setminus k) $ is 2-generated, we can choose the orders of two of these portals (and then
        the rest are determined); the covering induces an orbifold structure on the knot complement, where the orbifold has a knotted singular locus in the shape of $ k $.
        This is visualised in Thurston's talk \textit{Knots to Narnia} \cite{narnia}, and more modern visualisation work has been done by Brakke \cite{polycut}
        and S\"ummerman \cite{summermann21}.
\end{enumerate}

The generalised setting is now the following (where $\ZZ_n=\ZZ/n\ZZ$ denotes the cyclic group of order $ n $):
\begin{definition}
  Let $ a,b \in \hat{\NN} := \NN \cup \{\infty\} $ with $ \max\{a,b\} \geq 3 $. Define $ \Gamma_\rho^{a,b} $ to be the subgroup of $ \PSL(2,\CC) $ generated by
  \begin{displaymath}
    X = \begin{bmatrix} \exp(\pi i/a) & 1 \\ 0 & \exp(-\pi i/a) \end{bmatrix}, Y_\rho = \begin{bmatrix} \exp(\pi i/b) & 0 \\ \rho & \exp(-\pi i/b) \end{bmatrix}.
  \end{displaymath}
  The \textit{Riley slice}, which we will denote by $ \mc{R}^{a,b} $, is the set of $ \rho \in \CC $ such that $ \Gamma^{a,b}_\rho $ is discrete, isomorphic
  to $ \ZZ_a*\ZZ_b$, and such that $ \Omega(\Gamma^{a,b}_\rho)/\Gamma^{a,b}_\rho $ is a sphere with four cone points -- two of order $a$ and two of order $b$. When $a$ or $b$ is equal to $\infty$ the cone points are punctures,  or geometrically cusps, in the quotient surface.
\end{definition}

\begin{figure}
  \centering
  \begin{subfigure}{.24\textwidth}
    \centering
    \includegraphics[width=\textwidth]{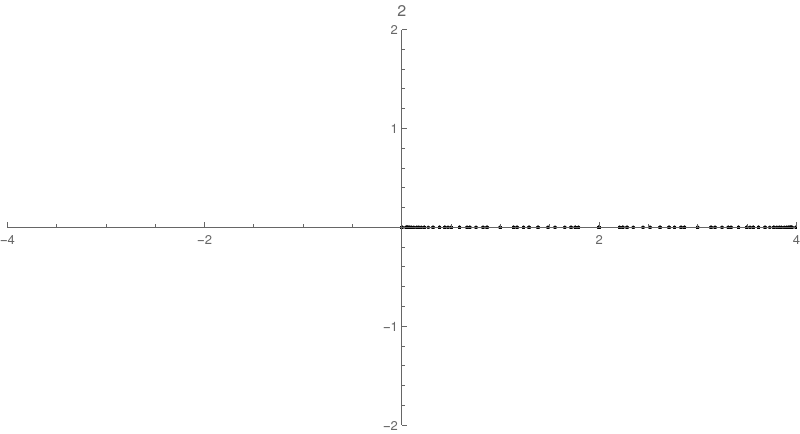}
    \caption{$(2,2)$}
  \end{subfigure}
  \begin{subfigure}{.24\textwidth}
    \centering
    \includegraphics[width=\textwidth]{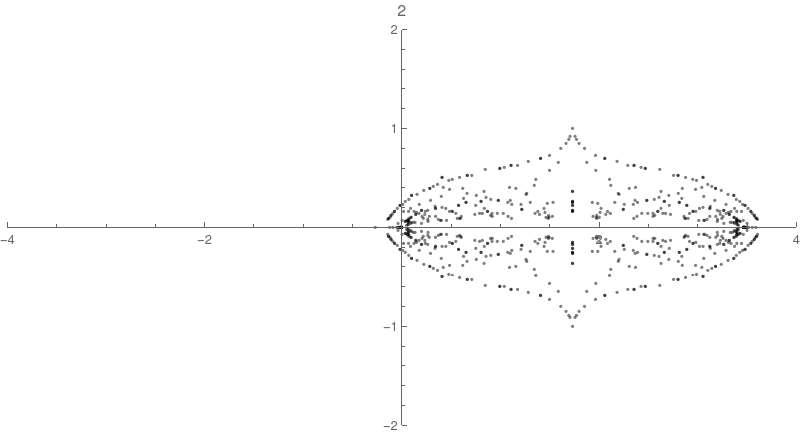}
    \caption{$(2,3)$}
  \end{subfigure}
  \begin{subfigure}{.24\textwidth}
    \centering
    \includegraphics[width=\textwidth]{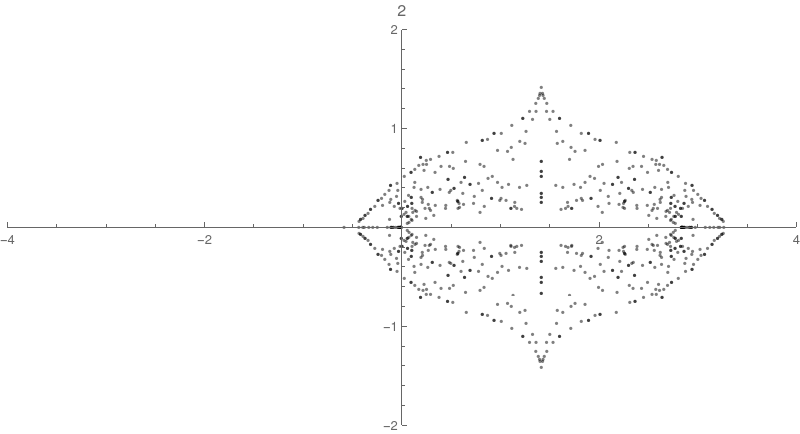}
    \caption{$(2,4)$}
  \end{subfigure}
  \begin{subfigure}{.24\textwidth}
    \centering
    \includegraphics[width=\textwidth]{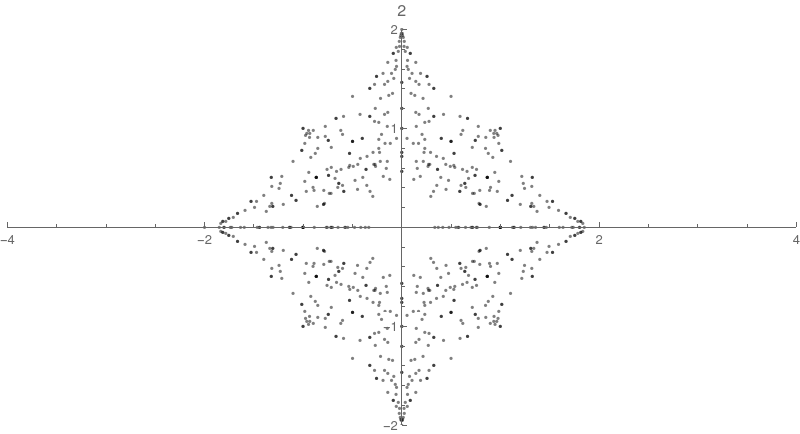}
    \caption{$(2,\infty)$}
  \end{subfigure}\\
  \begin{subfigure}{.24\textwidth}
    \centering
    \includegraphics[width=\textwidth]{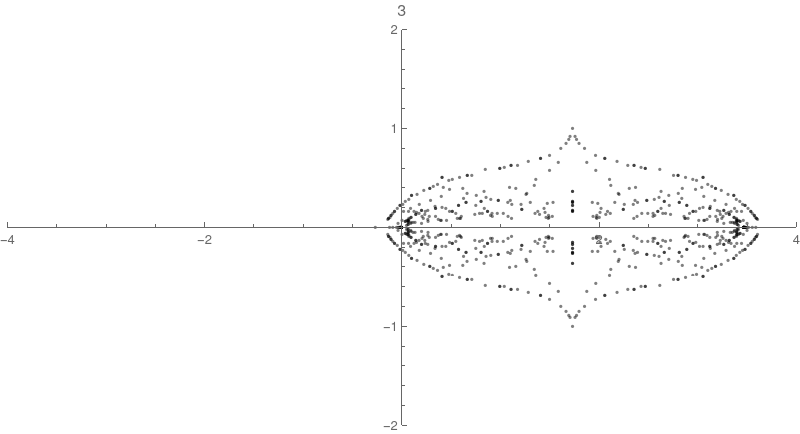}
    \caption{$(3,2)$}
  \end{subfigure}
  \begin{subfigure}{.24\textwidth}
    \centering
    \includegraphics[width=\textwidth]{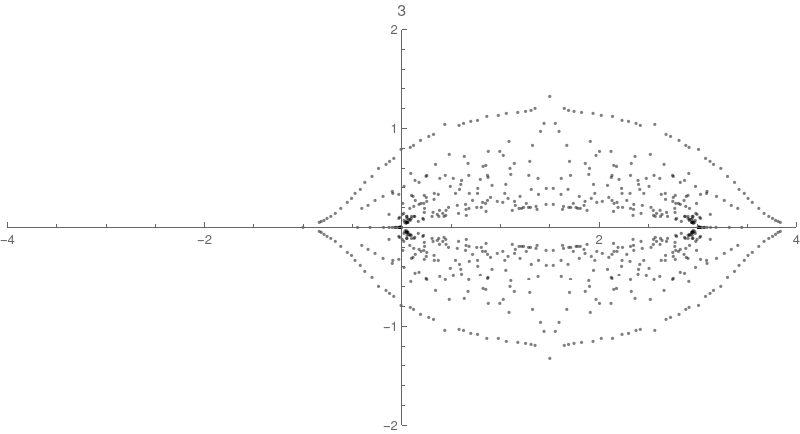}
    \caption{$(3,3)$}
  \end{subfigure}
  \begin{subfigure}{.24\textwidth}
    \centering
    \includegraphics[width=\textwidth]{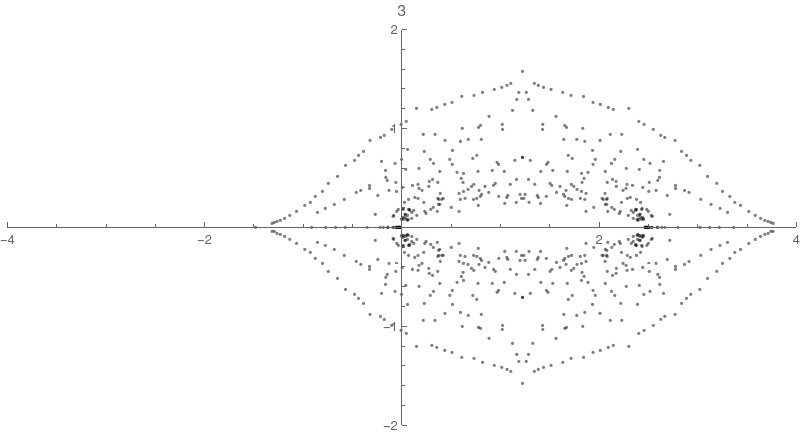}
    \caption{$(3,4)$}
  \end{subfigure}
  \begin{subfigure}{.24\textwidth}
    \centering
    \includegraphics[width=\textwidth]{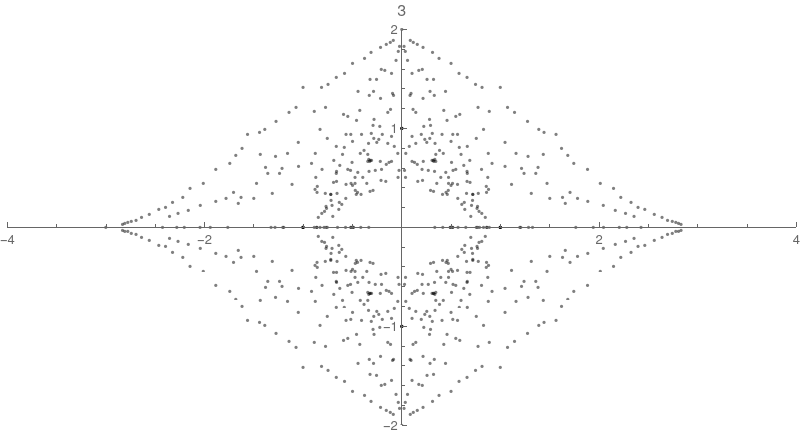}
    \caption{$(3,\infty)$}
  \end{subfigure}\\
  \begin{subfigure}{.24\textwidth}
    \centering
    \includegraphics[width=\textwidth]{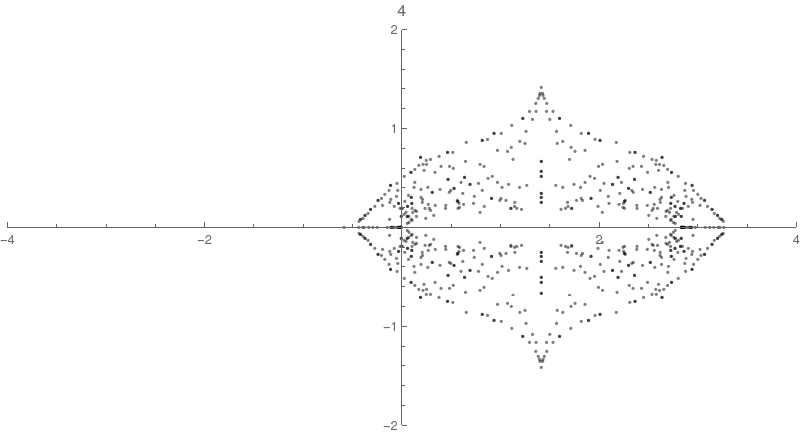}
    \caption{$(4,2)$}
  \end{subfigure}
  \begin{subfigure}{.24\textwidth}
    \centering
    \includegraphics[width=\textwidth]{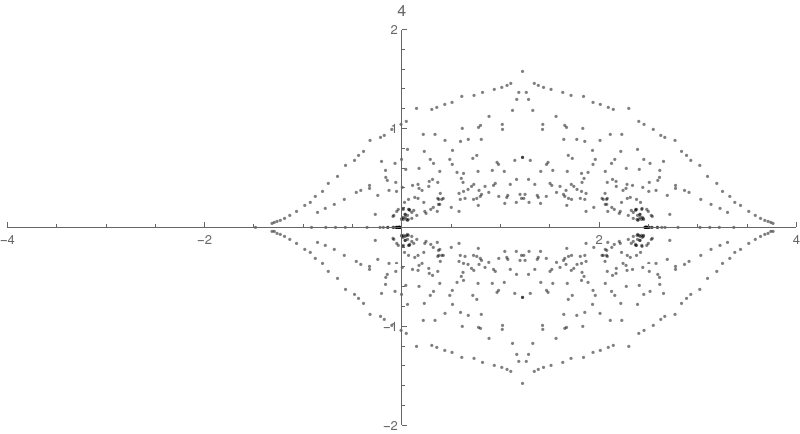}
    \caption{$(4,3)$}
  \end{subfigure}
  \begin{subfigure}{.24\textwidth}
    \centering
    \includegraphics[width=\textwidth]{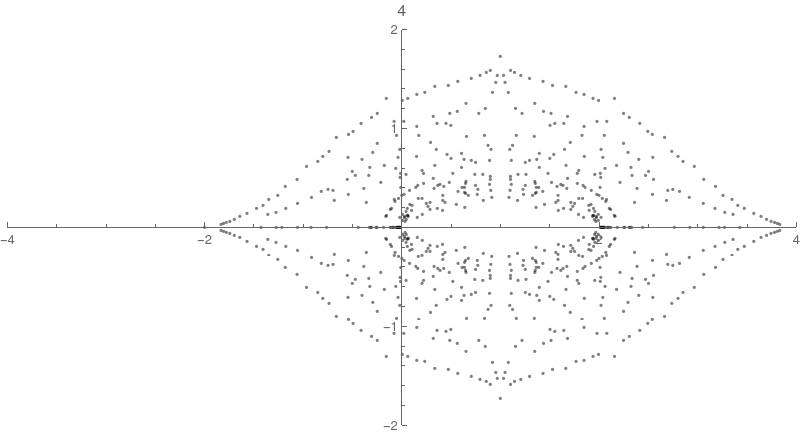}
    \caption{$(4,4)$}
  \end{subfigure}
  \begin{subfigure}{.24\textwidth}
    \centering
    \includegraphics[width=\textwidth]{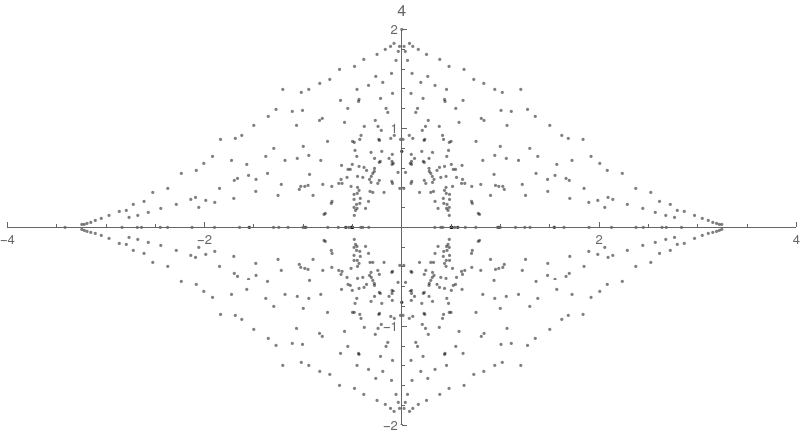}
    \caption{$(4,\infty)$}
  \end{subfigure}\\
  \begin{subfigure}{.24\textwidth}
    \centering
    \includegraphics[width=\textwidth]{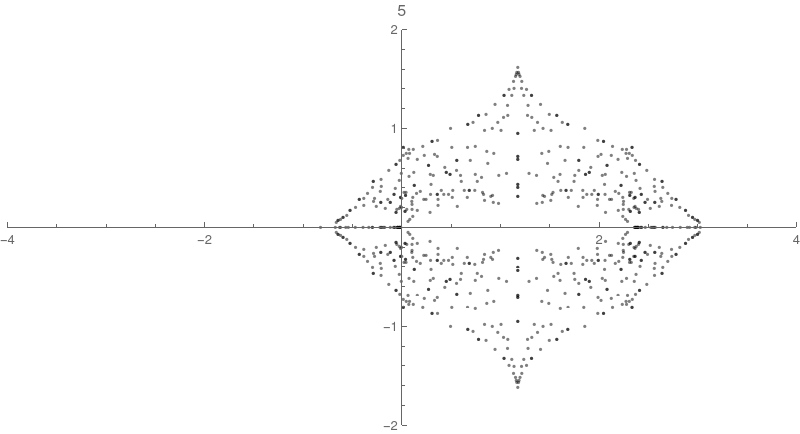}
    \caption{$(5,2)$}
  \end{subfigure}
  \begin{subfigure}{.24\textwidth}
    \centering
    \includegraphics[width=\textwidth]{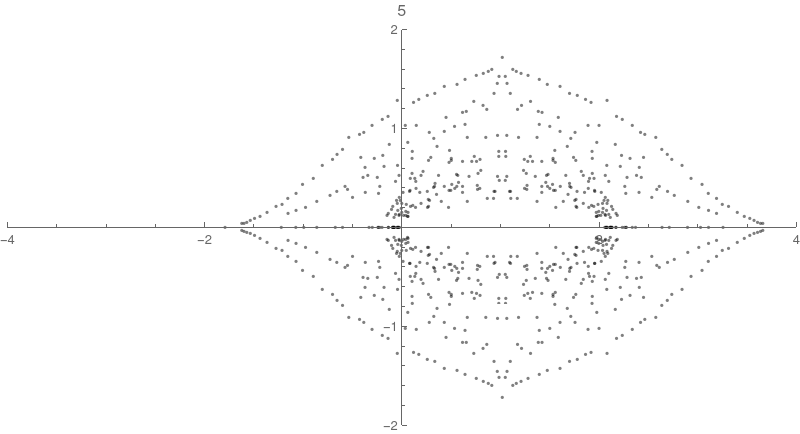}
    \caption{$(5,3)$}
  \end{subfigure}
  \begin{subfigure}{.24\textwidth}
    \centering
    \includegraphics[width=\textwidth]{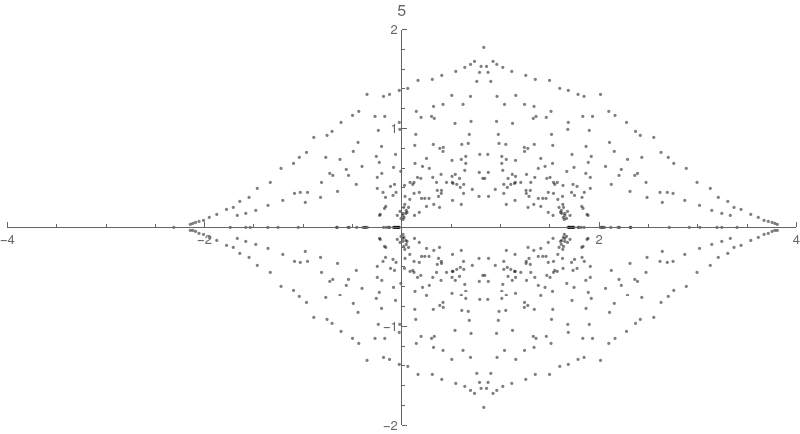}
    \caption{$(5,4)$}
  \end{subfigure}
  \begin{subfigure}{.24\textwidth}
    \centering
    \includegraphics[width=\textwidth]{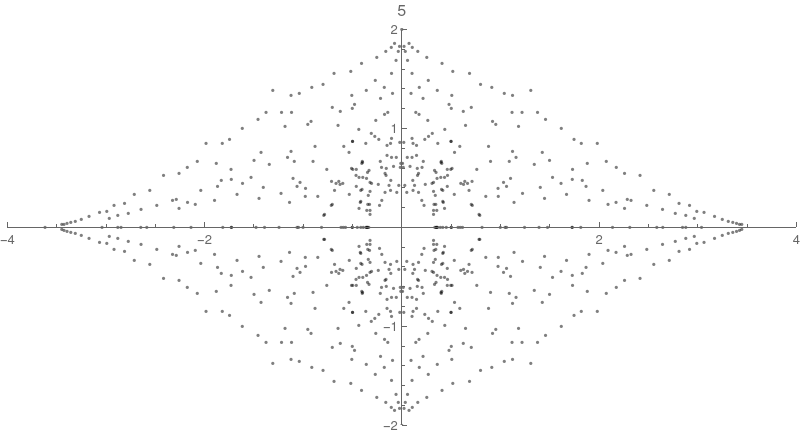}
    \caption{$(5,\infty)$}
  \end{subfigure}\\
  \begin{subfigure}{.24\textwidth}
    \centering
    \includegraphics[width=\textwidth]{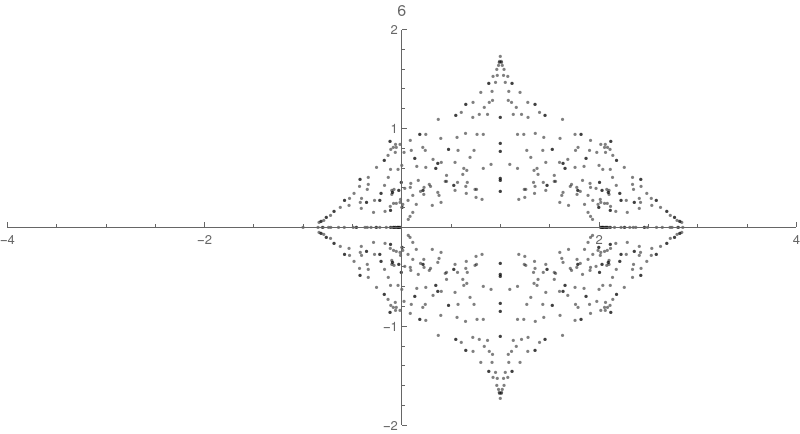}
    \caption{$(6,2)$}
  \end{subfigure}
  \begin{subfigure}{.24\textwidth}
    \centering
    \includegraphics[width=\textwidth]{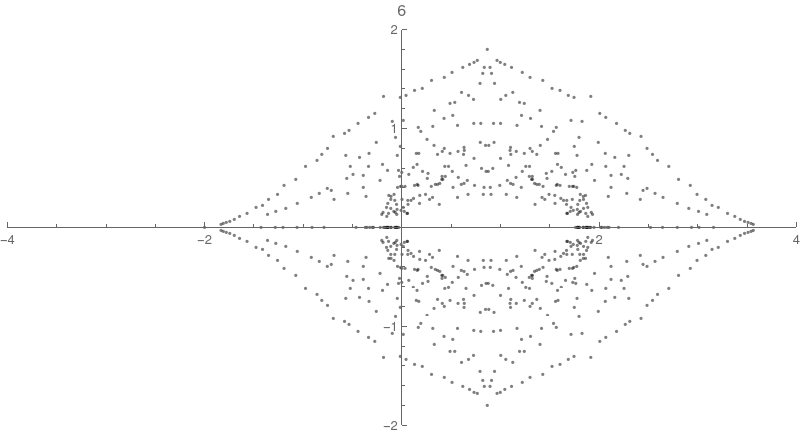}
    \caption{$(6,3)$}
  \end{subfigure}
  \begin{subfigure}{.24\textwidth}
    \centering
    \includegraphics[width=\textwidth]{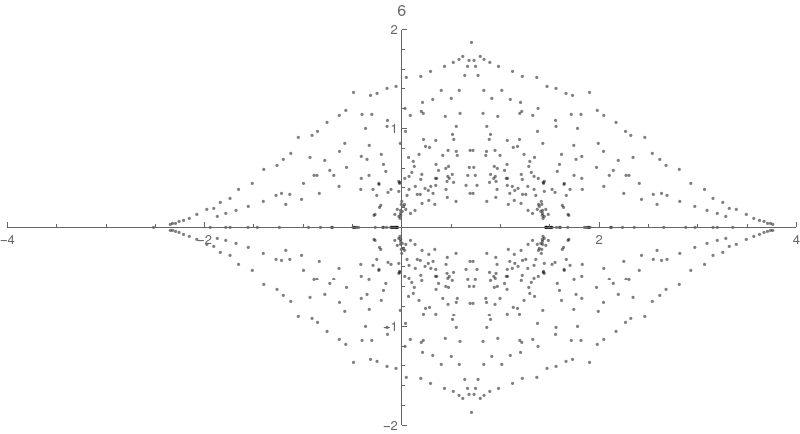}
    \caption{$(6,4)$}
  \end{subfigure}
  \begin{subfigure}{.24\textwidth}
    \centering
    \includegraphics[width=\textwidth]{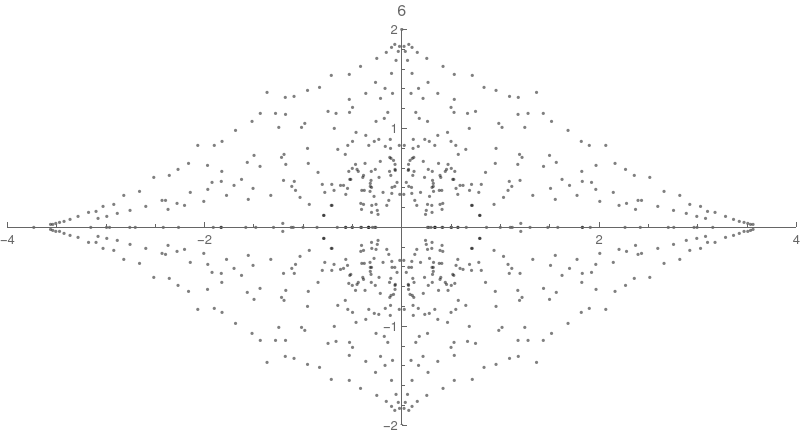}
    \caption{$(6,\infty)$}
  \end{subfigure}\\
  \begin{subfigure}{.24\textwidth}
    \centering
    \includegraphics[width=\textwidth]{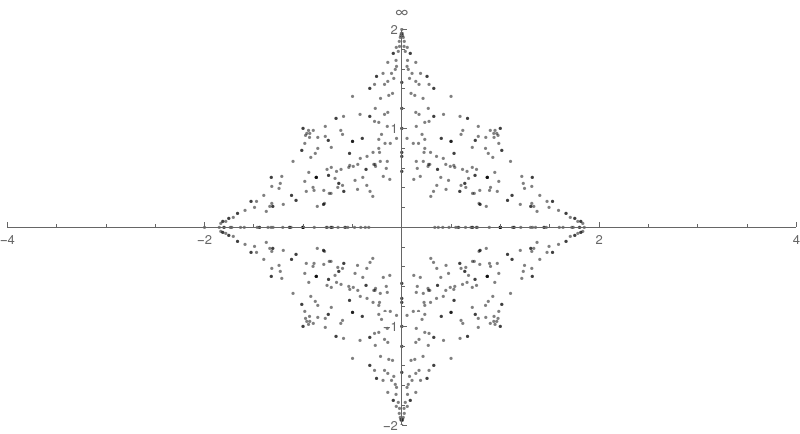}
    \caption{$(\infty,2)$}
  \end{subfigure}
  \begin{subfigure}{.24\textwidth}
    \centering
    \includegraphics[width=\textwidth]{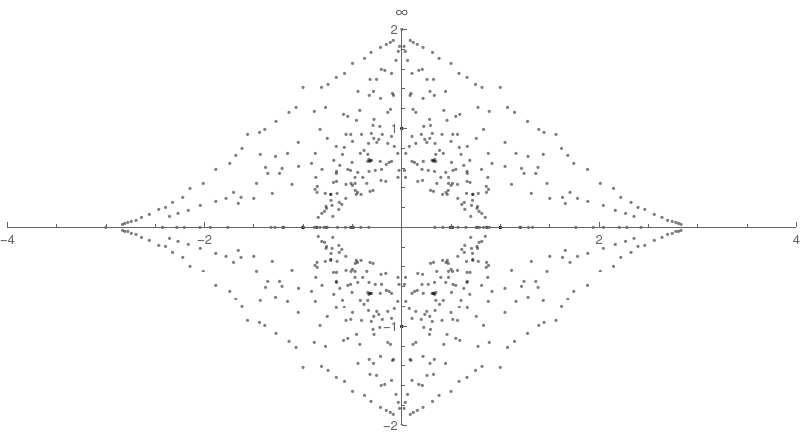}
    \caption{$(\infty,3)$}
  \end{subfigure}
  \begin{subfigure}{.24\textwidth}
    \centering
    \includegraphics[width=\textwidth]{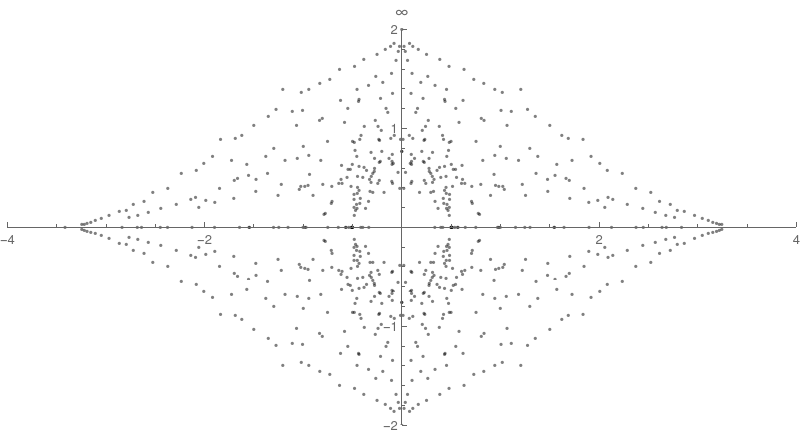}
    \caption{$(\infty,4)$}
  \end{subfigure}
  \begin{subfigure}{.24\textwidth}
    \centering
    \includegraphics[width=\textwidth]{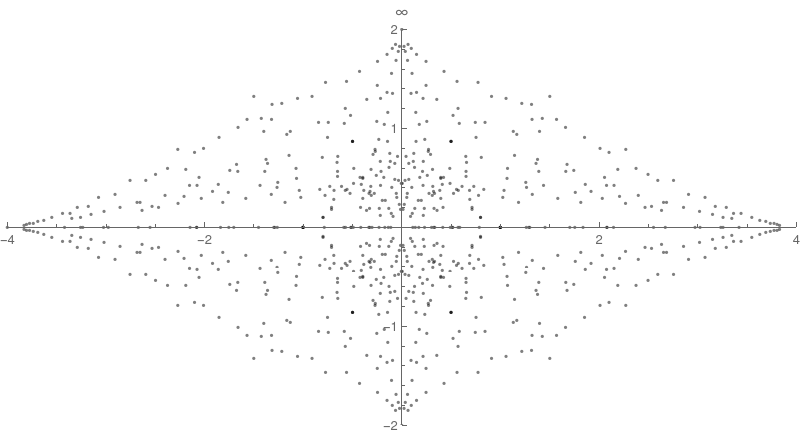}
    \caption{$(\infty,\infty)$}
  \end{subfigure}
  \caption{Riley slices for various cone point orders $ (a,b) $.\label{fig:riley_slice_elliptic}}   The $(2,2)$ picture is simply an interval as any group generated by two elements of order $2$ is virtually abelian and can be deformed until it becomes elliptic. 
\end{figure}

By convention, we take $ \pi i/\infty := 0 $; with this definition, $ \mc{R}^{\infty,\infty} $ is the classical parabolic Riley slice.
We draw pictures of various Riley slices (analogous to the parabolic slice drawn in Fig.~\ref{fig:crudebound} above) in Fig.~\ref{fig:riley_slice_elliptic}.

Only limited study of the elliptic case has been carried out to this point; for example, the PhD thesis \cite{zhang10}, and the study of
arithmetic groups in \cite{maclachlan99}. We remedy this in our preprint \cite{ems22}, where we justify why all of the topological theory
mentioned above in Section~\ref{sec:rileytop} carries through in the elliptic case, and carefully generalise the Keen--Series coordinate
system of Theorem~\ref{thm:ksmain} to the elliptic case. The main difficulty is in replacing topological arguments with orbi-analogues,
and generalising various hyperbolic arguments from ideal polyhedra to polyhedra with vertices contained in hyperbolic space. We also take
the opportunity to give more detail than Keen and Series gave in \cite{keen94} (who gave only sketches of some proofs which were similar
to those found in their earlier paper on the Maskit slice \cite{keen93}) and give more detailed motivation via knot theory in order to make
the theory more accessible to non-experts.

\medskip

Observing Figure \ref{fig:riley_slice_elliptic} one is struck by the apparent continuous deformations of these slices as the `variables' $a,b$ move.  One is tempted to conjecture that that there are natural quasiconformal deformations which transition between the various slices.  As these regions are defined dynamically or analytically (density of roots of Farey polynomials) this seems likely but we have not been able to establish this.

\subsection{Open neighbourhoods of pleating rays}\label{sec:nbds}
The theory of Keen and Series depends on deforming the limit sets of Kleinian groups in a conformal way. More precisely, and as they note, the rational pleating rays
are curves along which the combinatorics of certain patterns of circles in the limit set are preserved. Allowing small deformations off the rational
pleating rays also preserves these combinatorics; however, the circles now become quasicircles. In our preprint \cite{ems21}, we show that the large-scale geometry
of these quasicircles is controllable and that the resulting groups still lie in the Riley slice; thus we can construct open neighbourhoods of the pleating rays in
the Riley slice, using the Farey polynomials that are used to construct the rays themselves; these open neighbourhoods have well-behaved boundary, which gives some
control of the geometry of the Riley slice boundary. More precisely,  if $\alpha$ is a rational pleating ray in $ \mc{R}^{a,b} $,  then there is a rational slope $r/s$ and a Farey polynomial $p_{r/s}$ such that for a suitable choice of branch
\[ \alpha = p_{r/s}^{-1}((-\infty,-2]) \]
and $p^{-1}(-2)$ represents a cusp group in $ \partial \mc{R}^{a,b} $. With this same choice of branch we prove that
\[    p_{r/s}^{-1}(H) \subset \mc{R}^{a,b}, \quad H=\big\{x+i y: x \leq -2, \; y\in \mathbb{R}\setminus \{0\}\big\}  \]

\begin{figure} 
    \centering
    \includegraphics[width=\textwidth]{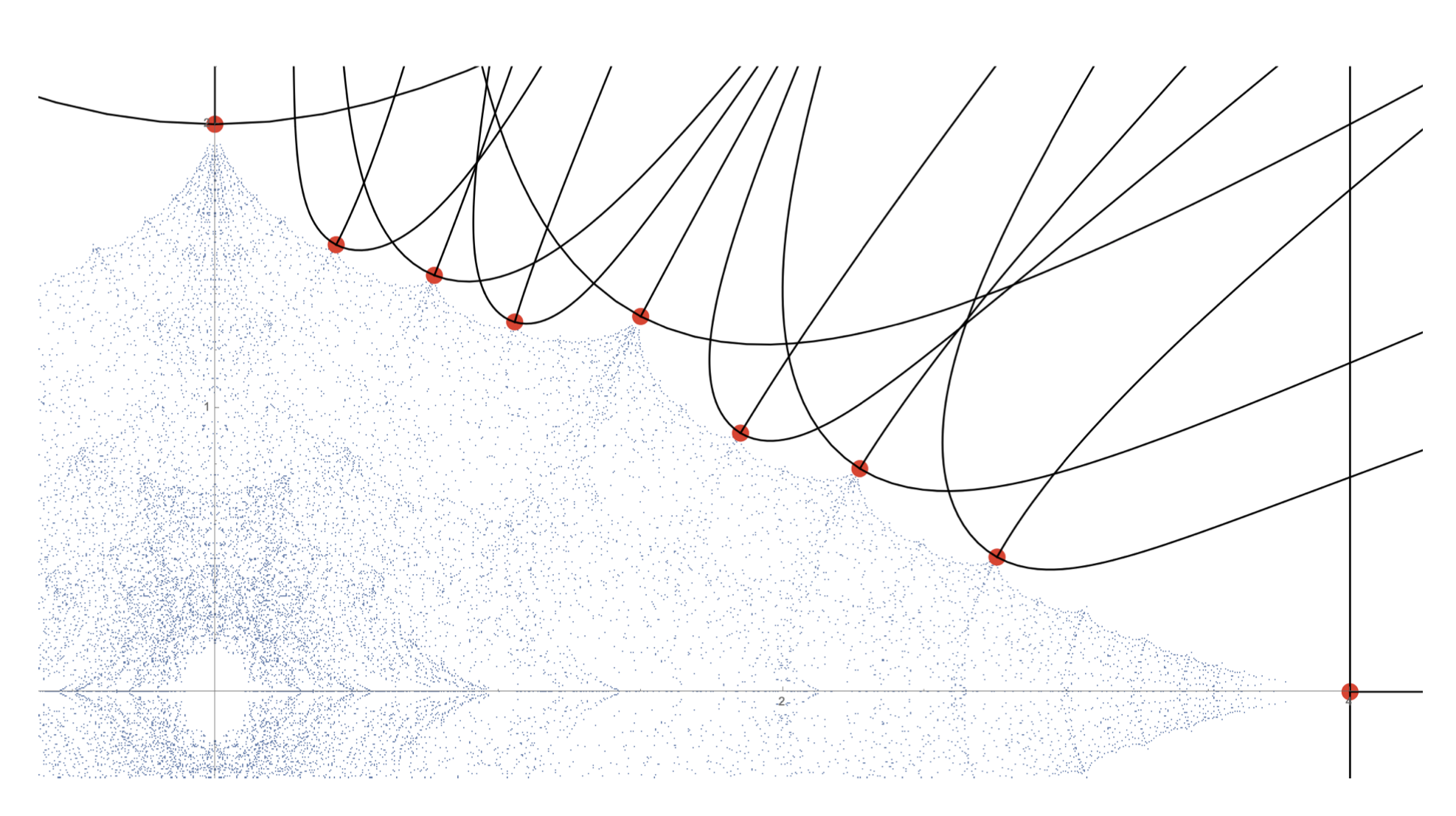}
  \caption{Neighbourhoods of pleating rays ``tangent'' to $\mc{R}^{a,b}$.}
\end{figure}

Because of the useful coarse description of $\mc{R}^{\infty,\infty}$ given by Lyndon and Ullman (see Fig.~\ref{fig:crudebound} and \cite{lyndon69}) a somewhat stronger result is true in that case \cite{ems21}.  We have been unable to completely replicate their result in the case of groups generated by two elliptics.

In the same preprint we also give some dynamical properties of the Farey polynomials, both for use
in studying these neighborhoods and for their own interest. This theory extends to the elliptic case in a similar way to the extension of the `round-circle'
Keen--Series theory, and the details will appear in the final section of \cite{ems22}.

\subsection{Farey polynomials}\label{sec:farey}
The Farey polynomials introduced by Keen and Series in \cite{keen94} also extend to the elliptic case, and we study their combinatorics in \cite{ems22b}. More precisely, we
derive a recursion for these polynomials (using results similar to those developed by Mumford, Series, and Wright \cite[p.277]{indras_pearls}, see also \cite{wright05}), show
formal analogies with the Chebyshev polynomials and other combinatorial objects, and give some possible applications to the approximation of the Riley slice boundary.

\bigskip

\begin{acknowledgement}
  The authors are grateful to the organisers of the MATRIX Workshop on Groups and Geometries in Nov--Dec 2021 and the subsequent workshop on Waiheke for the opportunity to speak on the Riley
  slice theory and our results. Some of the material in this article is based on sections of the first author's MSc thesis \cite{thesis}. Work of authors partially supported by the New Zealand
  Marsden Fund and by a University of Auckland FRDF grant.
\end{acknowledgement}

\end{document}